\documentclass[11pt,draft]{article}

\usepackage{amssymb,amsmath,amsthm}
\usepackage{epsfig}
\usepackage{verbatim}
\usepackage{graphicx}

\usepackage{color}

\newtheorem{thm}{Theorem}[section]

\newtheorem{rem}[thm]{Remark}

\newtheorem{lemma}[thm]{Lemma}

\newcommand{\R}{\mathbb{R}}

\newcommand{\D}{\displaystyle}

\newcommand{\esssup}{{\rm ess\thinspace sup}\thinspace}

\newcommand{\grad}{\nabla}

\newcommand{\F}{\mathcal{F}}

\newcommand{\dxu}{\partial_{x_1}}
\newcommand{\dxd}{\partial_{x_2}}

\newcommand{\dxi}{\partial_{x_i}}

\newcommand{\la}{\Lambda}

\newcommand{\sgn}{{\rm sgn}\thinspace}

\newcommand{\al}{\alpha}

\newcommand{\ep}{\varepsilon}

\newcommand{\Den}{\Delta_{-y}}

\newcommand{\eqdef}{\overset{\mbox{\tiny{def}}}{=}}

\begin{document}

\title{On the Muskat problem:\\ global in time results in 2D and 3D}

\author{Peter Constantin, Diego C{\'o}rdoba, Francisco Gancedo\\ Luis Rodr\'iguez-Piazza and Robert M. Strain}

\date{February 14, 2014}

\maketitle

\begin{abstract}
This paper considers the three dimensional Muskat problem in the stable regime. We obtain a conservation law which provides an $L^2$ maximum principle for the fluid interface. We also show global in time existence for strong and weak solutions with initial data controlled by explicit constants.  Furthermore we refine the estimates from our paper \cite{PDPB} to obtain global existence and uniqueness for strong solutions with larger initial data than we previously had in 2D. Finally we provide global in time results in critical spaces, giving solutions with bounded slope and time integrable bounded curvature.
\end{abstract}



\setcounter{tocdepth}{1}
%\tableofcontents

\section{Introduction}

   We consider the dynamics of the interface in between two incompressible fluids in porous media in three dimensional space. This is the Muskat problem (see \cite{Muskat}). We assume that both fluids are immiscible and
   have the same constant viscosity but different constant densities. We simplify matters by taking gravity $g = 1$, the permeability of the medium $\kappa = 1$ and the viscosity $\nu = 1$. Then the motion of the fluids satisfy:
   \begin{align}
\begin{split}\label{3DIPM}
\rho_t + \nabla\cdot(u\rho) &=0\\
u + \nabla P&=- (0,0,\rho)\\
\nabla\cdot u &=0
\end{split}
\end{align}
where $\rho=\rho(x_1,x_2,x_3,t)$ is the density, $P=P(x_1,x_2,x_3,t)$ is the pressure, $$u=(u_1(x_1,x_2,x_3,t),u_2(x_1,x_2,x_3,t), u_3(x_1,x_2,x_3,t))$$ is the incompressible velocity field, $x_i\in\R$, $i=$1, 2, 3 and $t\geq 0$. The first equation is the conservation of mass and the second equation is Darcy's law, where the velocity is proportional to the driving forces, the pressure gradient and the buoyancy force. We denote the interface that separates the space in two domains $\Omega_1$ and $\Omega_2$ by $x_3 = f(x_1,x_2,t)$. We consider the density $\rho=\rho(x_1,x_2,x_3,t)$ to be the following  step function:
\begin{equation}\label{density}
\rho(x_1,x_2,x_3,t)=\left\{\begin{array}{cl}
                    \rho^1,& (x_1,x_2,x_3)\in\Omega^1(t)=\{x_3> f(x_1,x_2,t)\}, \\
                    \rho^2  ,& (x_1,x_2,x_3)\in\Omega^2(t)=\{x_3< f(x_1,x_2,t)\}.
                 \end{array}\right.
\end{equation}
 Then the interface satisfies the equation
  \begin{align}
\begin{split}\label{ec}
\D f_t(x,t)&=\frac{\rho_2-\rho_1}{4\pi}PV\int_{\R^2}\frac{(\grad f(x,t)-\grad
f(x-y,t))\cdot y}{[|y|^2+(f(x,t)-f(x-y,t))^2]^{3/2}}dy,\\
f(x,0)&=f_0(x), \quad x=(x_1,x_2)\in\R^2,
\end{split}
\end{align}
in order to be a solution of the system (\ref{3DIPM}) (see \cite{DP} for a detail derivation).

A first approach is to linearize equation (\ref{ec}) around the steady state $f_0(x)=0$, which yields
\begin{align}
\begin{split}\label{ecl}
&f_t=\frac{\rho_1-\rho_2}{2}\la f,\\
&f(x,0)=f_0(x),
\end{split}
\end{align}
where the operator $\la f$ is defined in Fourier variables by $\widehat{\la
f}(\xi)=|\xi|\widehat{f}(\xi)$. The case $\rho_1<\rho_2$ gives a stable regime, and for $\rho_1>\rho_2$ the system is unstable. Stability versus instability is determined by the normal component of the pressure gradient jump at the interface  having a distinguished sign. This is known as the Rayleigh-Taylor condition which  implies local existence in $H^s$ when the heavier fluid is below the lighter one, and ill-posedness in the unstable regime (see \cite{DP} for a proof of both statements). Earlier works on the well-posedness in Sobolev spaces for the 3D Muskat problem, where both fluids have also different viscosities, include \cite{Am2}, \cite{SCH}, \cite{Escher} and \cite{ADP2}.

Our goal is to prove global in time existence results for the stable regime. Our main concern is about the size of the initial data needed to reach this conclusion. Global existence for large slopes turns out to be false. There exist initial data that turn to the unstable regime; in finite time the interface becomes no longer a graph (see \cite{CCFGL}). Moreover, there exist smooth initial data in the stable regime that in finite time turn to the unstable regime and at a later time they are no longer $C^4$ (see \cite{CCFG}). In our previous work \cite{PDPB}, we studied the two dimensional Muskat equation; we showed global existence of Lipschitz continuous solutions for initial data that satisfy $\|f_0\|_{L^\infty}<\infty$ and $\|\partial_x f_0\|_{L^\infty}<1$. We also proved  global existence  for unique strong solutions if the initial data is smaller than a constant $c_0$;  $\| f\|_1 \le c_0$ where
\begin{equation*}
 \|f\|_s
 \eqdef
 \int d\xi ~ |\xi|^s |\hat{f}(\xi)|,\quad s\geq 0.
\end{equation*}
We have checked numerically that $c_0$ is not \emph{small}; it is greater than $\frac15$. Recently, in \cite{Gr}, global results are obtained in a confined domain for initial data satisfying \emph{smallness} conditions relating the amplitude, the slope and the depth. We also point out a new work \cite{Beck} where instant analyticity is proved for small initial data represented on the Fourier side by positive measures.

In this paper we show that in 3D it is possible to obtain similar global existence results but with different constants. First, in Section \ref{sec:L2} we prove the following identity for the evolution of the $L^2$ norm of the contour
\begin{align}
\begin{split}\label{niceidentity}
\|f\|^2_{L^2}(t)+&\frac{\rho}{\pi}\int_0^t\!\!\int_{\R^2}\!\int_{\R^2}\!
\Big(\frac1{|y|}-\frac{1}{[|y|^2\!+\!(f(x,s)\!-\!f(x\!-\!y,s))^2]^{1/2}}\Big)dxdyds
=\|f_0\|^2_{L^2},
\end{split}
\end{align}
where $\rho=(\rho^2-\rho^1)/2$. We further explain using this this formula that there is no parabolic behavior in the contour equation at the level of $f$. In Section \ref{sec:17} we prove global existence of unique $C([0,T];H^k(\R^2))$ solutions for $k\geq 3$ if initially  $f_0$ is controlled by $\|f_0\|_1 < k_0$ where $k_0\geq 1/5$ (see \eqref{conko} for the exact size of $k_0$).   We also use the calculations in Section \ref{sec:17} to improve the size of the initial data in our global existence and uniqueness theorem for smooth solutions in \cite{PDPB} for 2D. In Section \ref{sec:13} we show that if a strong solution has the property $\|\grad f_0\|_{L^\infty}<1/3$, then it will be preserved in time. In Section \ref{sec:init1} we prove global in time existence of Lipschitz continuous solutions in the stable case for initial data satisfying $\|f_0\|_{L^\infty}<\infty$ and $\|\grad f_0\|_{L^\infty}<1/3$. Finally in Section 6 we use the parabolicity of the problem to show global in time solutions in critical spaces with $\|f_0\|_1<k_0$, $\|f\|_1(t)\in L^{\infty}([0,T])$ and $\|f\|_2(t)\in L^1([0,T])$ for any $T>0$. This result gives in particular that $\|f\|_{C^1}(t)\in L^{\infty}([0,T])$ and $\|f\|_{C^2}(t)\in L^1([0,T])$.

\section{$L^2$ maximum principle}\label{sec:L2}

This section is devoted to the proof of the identity (\ref{niceidentity}).
%\begin{align}
%\begin{split}\label{niceidentity}
%\|f\|^2_{L^2}(t)+\frac{\rho^2-\rho^1}{2\pi}&\int_0^t\int_{\R^2}\int_{\R^2}
%\frac{1}{|y|} \Big(1-\frac{1}{[1+(\Delta_yf(x))^2]^{1/2}}\Big)dxdyds\\
%&=\|f_0\|_{L^2},
%\end{split}
%\end{align}
%where
%$$\Delta_y f(x)=\frac{f(x)-f(x-y)}{|y|}.$$

In order to simplify the exposition we take $(\rho^2-\rho^1)/(2\pi)=1$ and we write $f(x,t)=f(x)$ for a fixed $t$. Then, the contour equation \eqref{ec} is given by:
\begin{align*}
\begin{split}
\D f_t(x)&=PV\int_{\R^2}\grad_x \Big(\frac{\Delta_yf(x)}{[1+(\Delta_yf(x))^2]^{1/2}}\Big)\cdot\frac{y}{|y|^2}dy,
\end{split}
\end{align*}
where
\begin{equation}\label{incre}
\Delta_yf(x)=(f(x)-f(x-y))/|y|.
\end{equation}

Integration by parts allows us to observe that
\begin{align*}
\begin{split}
\frac{1}{2}\frac{d}{dt}\|f\|^2_{L^2}(t)&=
-PV\int_{\R^2}  dy \int_{\R^2} dx ~
\grad f(x) \cdot \frac{y}{|y|^2} \frac{\Delta_yf(x)}{[1+(\Delta_yf(x))^2]^{1/2}}
\\
&=
-PV\int_{\R^2}  dy  \int_{\R^2} dx ~
\grad f(x) \cdot \frac{x-y}{|x-y|^2} \frac{\Delta_{x-y}f(x)}{[1+(\Delta_{x-y}f(x))^2]^{1/2}}.
\end{split}
\end{align*}
Next we split this into two terms,
\begin{align*}
\begin{split}
\frac{1}{2}\frac{d}{dt}\|f\|^2_{L^2}(t)&=-\int_{\R^2}dy\int_{\R^2}\frac{dx}{|x-y|}\frac{(\Delta_{x-y}f(x))^2}{[1+
(\Delta_{x-y}f(x))^2]^{1/2}}\\
&\quad-PV\int_{\R^2}dy\int_{\R^2}dx\frac{x-y}{|x-y|}\cdot \grad_x
([1+(\Delta_{x-y}f(x))^2]^{1/2}-1)\\
&=I_1+I_2.
\end{split}
\end{align*}
With these computations, a further integration by parts provides
$$
I_2=\int_{\R^2}\int_{\R^2}\frac{1}{|x-y|}
([1+(\Delta_{x-y}f(x))^2]^{1/2}-1)dx dy,
$$
and this equality gives
$$
\frac{1}{2}\frac{d}{dt}\|f\|^2_{L^2}=-\int_{\R^2}\int_{\R^2}
\frac{1}{|x-y|} \Big(1-\frac{1}{[1+(\Delta_{x-y}f(x))^2]^{1/2}}\Big)dx dy.
$$
From above \eqref{niceidentity} follows easily.

Next we show  the bound
$$
J\eqdef\int_{\R^2}\int_{\R^2}
\frac{1}{|y|} \Big(1-\frac{1}{[1+(\Delta_{y}f(x))^2]^{1/2}}\Big)dx dy\leq 4\pi\sqrt{2}\|f\|_{L^1},
$$
which controls the integral $J$ with zero derivatives. This expresses the fact that identity \eqref{niceidentity} does not give a gain of regularity at the level of $f$. Besides the linearization \eqref{ecl}, the nonlinear structure of the equation does not yield a parabolic dissipation for large initial data.

In order to deal with $J$ we observe that
$$
J\leq \int_{\R^2}\int_{\R^2} \frac{1}{|y|}\Big(1-\frac{1}{[1+(|f(x)|+|f(x-y)|)^2|y|^{-2}]^{1/2}}\Big) dx dy.
$$
Using the function $H(z)=1-(1+z^2)^{-1/2}$ and the fact that $$H(|z_1|+|z_2|)\leq H(\sqrt{2}|z_1|)+H(\sqrt{2}|z_2|)$$ it is easy to get
$$
J\leq \int_{\R^2}\int_{\R^2}\frac{1}{|y|}\Big(1\!-\!\frac{1}{[1\!+\!2|f(x)|^2|y|^{-2}]^{1/2}}\!+\!1\!-\!\frac{1}{[1\!+\!2|f(x\!-\!y)|^2|y|^{-2}]^{1/2}}\Big)dx dy,
$$
and therefore
$$
J\leq 2\int_{\R^2}\int_{\R^2}\frac{1}{|y|}
\Big(1-\frac{1}{[1+2|f(x)|^2|y|^{-2}]^{1/2}}\Big)dx dy=K.
$$
By an easy change of variable one finds
$$
K=2\sqrt{2}\int_{\R^2}|f(x)|dx \int_{\R^2}(\frac{1}{|z|}-\frac{1}{\sqrt{|z|^2+1}}) dz,
$$
so that $K=4\pi\sqrt{2}\|f\|_{L^1}.$ This provides the desired bound.

\section{A global existence result for data less than $\frac15$}\label{sec:17}

In this section we give a global existence result for classical solutions of the Muskat contour equation. We consider the norm
\begin{equation}
 \|f\|_s
 \eqdef
 \int_{\R^2}d\xi ~ |\xi|^s |\hat{f}(\xi)|,
 \quad
 s \ge 1,
 \label{norm:s}
\end{equation}
which allows us to use Fourier techniques for small initial data. We prove the following theorem:

\begin{thm}\label{17:thm}
Suppose that initially $f_0 \in H^l(\R^2)$ for $l\geq 3$ and $\|f_0\|_1 < k_0$, where $k_0$ is a constant such that
\begin{equation}\label{conko}
\pi\sum_{n\geq 1} (2n+1)^{1+\delta}\frac{(2n+1)!}{(2^n n!)^2}k_0^{2n}\leq 1,
\end{equation}
for some $0<\delta<1$. Then there is a unique solution $f$ of \eqref{ec} with initial data $f_0$ that satisfies $f \in C([0,T];H^l(\R^2))$ for any $T>0$.\end{thm}

\begin{rem}
Computing the limit case $\delta=0$, so that
$$
\pi\sum_{n\geq 1} (2n+1)\frac{(2n+1)!}{(2^n n!)^2}k_0^{2n}< 1
$$
one finds
$0<k_0\lesssim 0.24874641998890142626.$ In particular, this holds if $k_0 \le 1/5$.
\end{rem}

\begin{rem}
Analogous estimations allow us to obtain a better size for $\|f_0\|_1$ than in \cite{PDPB} in order to have a global existence and uniqueness result in 2D (1D interface). In fact, if initially $f_0 \in H^l(\R)$ for $l\geq 2$ and $\|f_0\|_1 < c_0$, where $c_0$ is a constant such that
\begin{equation}\label{conko}
2\sum_{n\geq 1} (2n+1)^{1+\delta}c_0^{2n}\leq 1,
\end{equation}
for some $0<\delta<1/2$, then there exists a unique solution $f$ of the two dimensional Muskat contour equation with initial data $f_0$ that satisfies $f \in C([0,T];H^l(\R))$ for any $T>0$. In the limit case $\delta=0$ we find $$0<c_0\leq \sqrt{(4-\sqrt{13})/3}\approx 0.3626057200026914,$$ and the result is true if for example $\|f_0\|_1\leq 1/3$.

\end{rem}

The remainder of this section is devoted to the proof of Theorem \ref{17:thm}. We point out that the argument used in \cite{PDPB} does not work directly here. It is valid in 2D only. To overcome the difficulty for 3D we need to symmetrize the operators involved in the equation to find an extra cancellation.
We define $\Delta_y f(x)$ as in \eqref{incre}
and we take $(\rho^2-\rho^1)/2=1$ for the sake of simplicity.
The contour equation for the Muskat problem \eqref{ec} can be written as
\begin{equation}
f_t(x,t)=-\Lambda f -N(f),
\label{muskatEQ2d}
\end{equation}
where the operator $\Lambda$ is the square root of the negative Laplacian and we have
\begin{align*}
N(f)=\frac1{2\pi}\int_{\R^2}\frac{y}{|y|^2}\cdot \grad_x \Delta_y f(x)
R(\Delta_y f(x))dy,
\end{align*}
with $R(z)=1-1/(1+z^2)^{3/2}.$ A change of variable allows us to obtain
\begin{align}
\begin{split}\label{td}
N(f)=\frac1{4\pi}\int_{\R^2}\frac{y}{|y|^2}\cdot \Big(\grad_x \Delta_y f(x) & R(\Delta_y f(x))\\
&-\grad_x \Den f(x) R(\Den f(x))\Big)dy.
\end{split}
\end{align}
We consider the norm
$
\|f\|_1
$
\eqref{norm:s} as follows:
\begin{align*}
\frac{d}{dt}\|f\|_1(t)&=\int_{\R^2}d\xi ~ |\xi| ~ (\hat{f}_t(\xi)\overline{\hat{f}(\xi)}+\hat{f}(\xi)\overline{\hat{f}_t(\xi)})/(2|\hat{f}(\xi)|)
\\
&\leq -\int_{\R^2}d\xi ~ |\xi|^2|\hat{f}(\xi)|+\int_{\R^2}d\xi ~|\xi||\F(N(f))(\xi)|.
\end{align*}
We will show that the first term controls the evolution in such a way that $\|f\|_1$ is decreasing if initially $$\|f_0\|_1< k_0,\mbox{ where } k_0\approx 0.24874641998890142626.$$ Since $|\Delta_y f(x)|\leq  \|f\|_1<1$ we can use the Taylor expansion
$$
R(z)=-\sum_{n\geq1} (-1)^{n}a_n z^{2n},\quad\mbox{with}\quad a_n=\frac{(2n+1)!}{(2^n n!)^2},\quad |z|<1,
$$
to obtain
\begin{align}
\begin{split}\label{tdt}
N(f)=\frac{-1}{4\pi}\D\sum_{n\geq 1}(-1)^n & a_n\int_{\R^2} \frac{y}{|y|^2}\cdot \\
\times & \Big(\grad_x (\Delta_y f)\,  (\Delta_y f)^{2n}-\grad_x (\Den f)\,(\Den f)^{2n} \Big)dy.
\end{split}
\end{align}
Recall that
\begin{gather}\notag
\F(\Delta_y f)=\hat{f}(\xi) m(\xi,y),\qquad \F(\grad_x \Delta_y f)=i\xi \hat{f}(\xi) m(\xi,y),
\\
\mbox{where}\qquad m(\xi,y)=(1-e^{-i\xi\cdot y})/|y|.
\end{gather}
Therefore
$$
\F(\grad_x (\Delta_y f)\,  (\Delta_y f)^{2n})=((i\xi\hat{f} m)\ast (\hat{f} m)\ast \cdots \ast (\hat{f} m))(\xi,\al),
$$
with $2n$ convolutions, one with $i\xi\hat{f} m$ and $2n-1$ with $\hat{f} m$.
Using \eqref{tdt}
\begin{align*}
\F(N)(\xi)
=&
\frac{-i}{4\pi}\D\sum_{n\geq 1}(-1)^na_n\int_{\R^2} dy \int_{\R^2} d\xi_1\cdots\int_{\R^2} d\xi_{2n} \frac{y}{|y|^2}\cdot (\xi\!-\!\xi_1)
\\
&\times\hat{f}(\xi\!-\!\xi_{1})\Big(\prod_{j=1}^{2n-1}\hat{f}(\xi_j\!-\!\xi_{j+1})\Big)
\hat{f}(\xi_{2n}) (M_n(y)-M_n(-y)),
\end{align*}
where $M_n(y)=M_n(\xi,\xi_1,\ldots,\xi_{2n},y)$ is given by
$$
M_n(y)=m(\xi\!-\!\xi_1,y)\Big(\prod_{j=1}^{2n-1}m(\xi_j\!-\!\xi_{j+1},y)\Big)m(\xi_{2n},y).
$$
We then use Fubini theorem to obtain
\begin{align}
\begin{split}\label{FN}
\F(N)(\xi)=&\D\sum_{n\geq 1}a_n
\int_{\R^2}d\xi_1\cdots\int_{\R^2} d\xi_{2n}\\
&\times(\xi\!-\!\xi_1)\hat{f}(\xi\!-\!\xi_{1})\Big(
\prod_{j=1}^{2n-1}  \hat{f}(\xi_j -\!\xi_{j+1})\Big)\hat{f}(\xi_{2n}) \cdot  I_n,
\end{split}
\end{align}
where the integral $I_n=I_n(\xi,\xi_1,\ldots,\xi_{2n})$ reads
\begin{align*}
I_n\eqdef\frac{-i}{4\pi}(-1)^n\!\int_{\R^2}\!\frac{y}{|y|^2} (M_n(y)-M_n(-y)) dy.
\end{align*}
Polar coordinates, $y=ru$ with $u=(\cos \theta,\sin \theta)$, provide
\begin{align*}
I_n=\frac{-i}{4\pi}(-1)^n\int_{-\pi}^{\pi}ud\theta\int_0^{+\infty}dr  (M_n(r,u)-M_n(r,-u)),
\end{align*}
where we redefine
$$
m(\xi,r,u)=(1-e^{-ir\xi\cdot u})/r,
$$
and
$$
M_n(r,u)=m(\xi\!-\!\xi_1,r,u)\Big(\prod_{j=1}^{2n-1}m(\xi_j\!-\!\xi_{j+1},r,u)\Big)m(\xi_{2n},r,u).
$$
Since we have $m(\xi,-r,u)=-m(\xi,r,-u)$ and $-m(\xi,-r,-u)=m(\xi,r,u)$, the change of variable $r=-s$ yields
\begin{align*}
I_n=\frac{-i}{4\pi}(-1)^n\int_{-\pi}^{\pi}u\, d\theta \int_{-\infty}^0 ds (M_n(s,u)-M_n(s,-u)),
\end{align*}
and therefore
\begin{align*}
I_n=\frac{-i}{8\pi}(-1)^n\int_{-\pi}^{\pi}u\,d\theta\int_\R dr  (M_n(r,u)-M_n(r,-u)).
\end{align*}
The identity $m(\xi,r,u) = i\xi\cdot u \int_0^1 ds~ e^{ir (s-1)\xi\cdot u }$
allows us to obtain
\begin{align*}
M_n(r,u)&=(-1)^n\int_0^1\!\!ds_1\cdots\int_0^1\!\!ds_{2n}\,
\left(\prod_{j=1}^{2n-1}(\xi_j-\xi_{j+1})\cdot u  \right)
\, \xi_{2n}\cdot u\\
&\times\frac{1\!-\!e^{-ir(\xi-\xi_1)\cdot u }}{r}\exp\Big(ir\Big(\sum_{j=1}^{2n-1}(s_j\!-\!1)(\xi_j\!-\!\xi_{j+1})\!+\!(s_{2n}\!-\!1)\xi_{2n}\Big)\!\cdot\! u\Big),
\end{align*}
which is simplified by writing
\begin{align*}
M_n(r,u)&=(-1)^n\int_0^1\!\!ds_1\cdots\int_0^1\!\!ds_{2n}\,
\left(\prod_{j=1}^{2n-1}(\xi_j-\xi_{j+1})\cdot u  \right)  \, \xi_{2n}\cdot u\\
&\times \Big(\frac{\exp(ir A\cdot u)}{r}-\frac{\exp(ir B\cdot u)}{r}\Big),
\end{align*}
with
\begin{align*}
A=\sum_{j=1}^{2n-1}(s_j-1)(\xi_j-\xi_{j+1})+(s_{2n}-1)\xi_{2n},
\end{align*}
and
\begin{align*}
B= -(\xi-\xi_1)+\sum_{j=1}^{2n-1}(s_j-1)(\xi_j-\xi_{j+1})+(s_{2n}-1)\xi_{2n}.
\end{align*}
It follows that
\begin{align*}
I_n&=\frac{-i}{8\pi}\int_{-\pi}^{\pi}ud\,\theta\int_0^1\!\!ds_1\cdots\int_0^1\!\!ds_{2n}\,
\left(\prod_{j=1}^{2n-1}(\xi_j-\xi_{j+1})\cdot u  \right)
\, \xi_{2n}\cdot u\\
\times& \int_\R dr  \Big(\frac{\exp(ir A\cdot u)}{r}-\frac{\exp(ir B\cdot u)}{r}-\frac{\exp(-ir A\cdot u)}{r}+\frac{\exp(-ir B\cdot u)}{r}\Big)
\end{align*}
and the equality $PV\int_\R dr \exp(ir\al)/r=\pi i \sgn \al$ yields
\begin{align*}
I_n=\frac{-1}{4}\int_{-\pi}^{\pi}ud\,\theta\int_0^1\!\!ds_1\cdots\int_0^1\!\!ds_{2n}\,
\Big(&\prod_{j=1}^{2n-1}(\xi_j  -\xi_{j+1})\cdot u  \Big)\\
\times  &\xi_{2n}\cdot u(\sgn (A\cdot u)-\sgn (B\cdot u)).
\end{align*}
At this point it is easy to bound $I_n$:
$$
|I_n|\leq \pi \prod_{j=1}^{2n-1}|\xi_j-\xi_{j+1}||\xi_{2n}|.
$$
The above estimate and \eqref{FN} allow us to get
\begin{align*}
\int_{\R^2}d\xi ~ |\xi|&|\F(N)(\xi)|\leq \pi\D\sum_{n\geq 1}a_n\int_{\R^2}d\xi\int_{\R^2} d\xi_1\cdots\int_{\R^2} d\xi_{2n} ~
|\xi|\\
\times&  |\xi-\xi_1||\hat{f}(\xi-\xi_{1})| \prod_{j=1}^{2n-1}  |\xi_j -\!\xi_{j+1}||\hat{f}(\xi_j -\!\xi_{j+1})||\xi_{2n}||\hat{f}(\xi_{2n})|.
\end{align*}
The inequality $|\xi|\leq |\xi-\xi_1|+|\xi_1-\xi_2|+\cdots+|\xi_{2n-1}-\xi_{2n}|+|\xi_{2n}|$ yields
\begin{align*}
\int_{\R^2}d\xi |\xi||\F(N)(\xi)| &\leq \pi\sum_{n\geq
1}(2n+1)a_n\Big(\int_{\R^2}d\xi |\xi|^2 |\hat{f}(\xi)|\Big)\Big(\int_{\R^2}d\xi |\xi| |\hat{f}(\xi)|\Big)^{2n} ,
\end{align*}
and therefore
\begin{align*}
\int_{\R^2}d\xi |\xi||\F(N)(\xi)| &\leq  \Big(\int_{\R^2}d\xi |\xi|^2 |\hat{f}(\xi)|\Big)\pi\sum_{n\geq
1}(2n+1)a_n\|f\|_1^{2n}\\
&\leq \Big(\int_{\R^2}d\xi |\xi|^2 |\hat{f}(\xi)|\Big) \pi\Big(\frac{1+2\|f\|^2_1}{(1-\|f\|^2_1)^{5/2}}-1\Big).
\end{align*}
For $0\leq x< k_0\approx 0.2487461998890142626$ one finds $(1+2x^2)/(1-x^2)^{5/2}-1<1/\pi$. Therefore if $\|f_0\|_1<k_0$ this inequality will be maintained when we propagate forward in time because of
$$
\frac{d}{dt}\|f\|_1(t)\leq 0,
$$
and $\|f\|_1(t) \le \|f_0\|_1<k_0$.

Considering a higher order norm, with $s>1$ in \eqref{norm:s}, we aim to obtain
\begin{equation}\label{qsmp}
\|f\|_{1+\delta}(t)+\mu\int_0^t ds ~ \|f\|_{2+\delta}(s)\leq \|f_0\|_{1+\delta},
\end{equation}
for some $0<\delta<1$ and $0<\mu<1$. Let us recall that
$
\|f_0\|_{1+\delta}\leq C\|f_0\|_{H^3}$
for $0<\delta<1$. We use the inequality
$$
|\xi|^{1+\delta}\leq (2n+1)^{\delta}(|\xi-\xi_1|^{1+\delta}+|\xi_1-\xi_2|^{1+\delta}+\cdots+
|\xi_{2n-1}-\xi_{2n}|^{1+\delta}+|\xi_{2n}|^{1+\delta}),
$$
to obtain as before
$$
\int_{\R^2}|\xi|^{1+\delta}|\F(N)(\xi)|d\xi\leq \int_{\R^2}|\xi|^{2+\delta}
|\hat{f}(\xi)| d\xi \, \pi\sum_{n\geq 1}(2n+1)^{1+\delta}a_n\|f\|_1^{2n}.
$$
Due to
$$
1>\pi\sum_{n\geq 1} (2n+1)^{1+\delta} a_n\|f_0\|_1^{2n}=1-\mu\geq \pi\sum_{n\geq 1} (2n+1)^{1+\delta} a_n\|f\|_1^{2n}(t),
$$
for some $0<\mu<1$, we find
$$
\int_{\R^2}|\xi|^{1+\delta}|\F(N)(\xi)|d\xi\leq (1-\mu) \int_{\R^2}|\xi|^{2+\delta}
|\hat{f}(\xi)| d\xi,
$$
for $\delta$ small enough. Since
$$
\frac{d}{dt}\|f\|_{1+\delta}(t)\leq -\mu \|f\|_{2+\delta}(t),
$$
integration in time provides \eqref{qsmp}.

From previous work \cite{DP}, one could find the following a priori bound:
\begin{equation*}
\frac{1}{2}\frac{d}{dt}\|\partial^3_{x_1}f\|^2_{L^2}\leq P(\|\grad f\|_{L^\infty})(\|\grad^2 f\|_{L^\infty}|\grad f|_{C^\delta}+\|\grad f\|_{L^\infty}|\grad^2 f|_{C^\delta})\|f\|^2_{H^3},
\end{equation*}
where $P$ is a polynomial function and $|\cdot|_{C^\delta}$ is the homogeneous H\"{o}lder norm. The terms that appear in the evolution can be handled as in \cite{DP} (see Section 4) except for a couple of low order terms:
$$
L.O.T.^1=\int_{\R^2}\partial^3_{x_1}f(x)\int_{\R^2}\frac{\grad_x\Delta_yf(x)\cdot y}{|y|^2}\frac{(\Delta_yf(x))^3(\partial_{x_1}\Delta_yf(x)))^3}{[1+(\Delta_yf(x))^2]^{9/2}}dydx,
$$
$$
L.O.T.^2=\int_{\R^2}\partial^3_{x_1}f(x)\int_{\R^2}\frac{\grad_x\Delta_yf(x)\cdot y}{|y|^2}\frac{(\Delta_yf(x))(\partial_{x_1}\Delta_yf(x)))^3}{[1+(\Delta_yf(x))^2]^{7/2}}dydx.
$$
We bound
$$
L.O.T.^1+L.O.T.^2\leq 2\int_{\R^2}|\partial^3_{x_1}f(x)|\int_{\R^2}\frac{|\grad_x\Delta_yf(x)|}{|y|}|\partial_{x_1}\Delta_yf(x))|^3dydx=J.
$$
Splitting $J$ for $|y|>1$ and $|y|<1$ it is easy to find
\begin{align*}
J&= \int_{\R^2}dx\int_{|y|>1}dy+\int_{\R^2}dx\int_{|y|<1}dy\\
&\leq C\|\partial^3_{x_1}f\|_{L^2} \|\grad f\|_{L^2}
\|\grad f\|_{L^\infty}|\grad f|_{C^\delta}\| \grad^2 f\|_{L^\infty}\int_{|y|>1}\frac{|y|^\delta}{|y|^4}dy\\
&\quad+C\|\partial^3_{x_1}f\|_{L^2} \|\grad^2 f\|^2_{L^4}
|\grad f|_{C^\delta}\|\grad^2 f\|_{L^\infty}\int_{|y|<1}\frac{|y|^\delta}{|y|^2}dy.
\end{align*}
Interpolation inequality $\|\grad^2 f\|^2_{L^4}\leq
\|\grad f\|_{L^\infty}\|\grad^3 f\|_{L^2}$ allows us to obtain
$$
J\leq C\|\grad f\|_{L^\infty}|\grad f|_{C^\delta}\|\grad^2 f\|_{L^\infty}\|f\|^2_{H^2},
$$
as desired. Proceeding in a similar way for $\|\partial^3_{x_2}f\|_{L^2}$ we find
$$
\frac{d}{dt}\|f\|^2_{H^3}\leq P(\|\grad f\|_{L^\infty})(\|\grad^2 f\|_{L^\infty}|\grad f|_{C^\delta}+\|\grad f\|_{L^\infty}|\grad^2 f|_{C^\delta})\|f\|^2_{H^3}.
$$
Fourier transform yields $\|\grad^k f\|_{L^\infty}\leq \|f\|_{k}$ and $|\grad^k f|_{C^\delta}(t)\leq\|f\|_{k+\delta}$ for $k=1,2$ and by interpolation it is easy to obtain
\begin{equation*}
\|f\|_{2}\|f\|_{1+\delta}\leq \|f\|_{1}\|f\|_{2+\delta}.
\end{equation*}
We find
$$
\frac{d}{dt}\|f\|^2_{H^3}\leq P(\|f\|_{1}) \|f\|_{2+\delta}\|f\|^2_{H^3},
$$
which together with the a priori bound provides
$$
\|f\|_{H^3}(t)\leq \|f_0\|_{H^3}\exp (CP(k_0)\int_0^t\|f\|_{2+\delta}(s)ds),
$$
after integration in time. Using \eqref{qsmp} we get finally
\begin{equation*}\label{h3n}
\|f\|_{H^3}(t)\leq \|f_0\|_{H^3}\exp (CP(k_0)\|f_0\|_{1+\delta}/\mu).
\end{equation*}
We finish with the conclusion that the solution can be continued in $H^3$ for all time if $\|f_0\|_1$ is initially smaller than $k_0$ defined by \eqref{conko}. An analogous calculation gives
$$
\|f\|_{H^k}(t)\leq \|f_0\|_{H^k}\exp (CP(k_0)\|f_0\|_{1+\delta}/\mu),
$$
getting the result for any $H^k$ for $k> 3$.

\section{Initial data smaller than 1/3}\label{sec:13}

In this section our goal is to prove the following maximum principle for the evolution of $\|\grad f\|_{L^\infty}(t)$ assuming that $\|\grad f_0\|_{L^\infty}<1/3$.

\begin{thm}
Let $f_0\in H^s$ with $s\geq 4$ and $\|\grad f_0\|_{L^\infty}<1/3$. Then the unique solution of the system \eqref{ec} satisfies
$$
\|\grad f\|_{L^\infty}(t)<1/3, \quad \mbox{for} \quad t>0.
$$
\end{thm}

Proof: We consider $(\rho^2-\rho^1)/2=1$ without loss of generality. We take one derivative in $x_i$ in \eqref{ec} to find
$$
\dxi f_t(x,t)=I^i_1(x,t)+I^i_2(x,t)+I^i_3(x,t),
$$
where
$$
I^i_1=\frac{1}{2\pi}PV\int_{\R^2} \frac{\grad \dxi f(x,t)\cdot y}{[|y|^2+(f(x,t)-f(x-y,t))^2]^{3/2}}dy,
$$
$$
I^i_2=-\frac{1}{2\pi}PV\int_{\R^2}\frac{\grad
\dxi f(x-y,t)\cdot y}{[|y|^2+(f(x,t)-f(x-y,t))^2]^{3/2}}dy
$$
and
$$
I^i_3=-\frac{1}{2\pi}\int_{\R^2}\frac{\dxi f(x,t)\!-\!\dxi f(x\!-\!y,t)}{[|y|^2\!+\!(f(x,t)\!-\!f(x\!-\!y,t))^2]^{3/2}}A(x,y)dy,
$$
with
$$
A(x,y)=3\frac{(f(x,t)\!-\!f(x-y,t))(\grad f(x,t)-\grad f(x-y,t))\cdot y}{|y|^2+(f(x)-f(x-y))^2}.
$$
Integration by parts yields
\begin{align*}
I^i_2=&-\frac{1}{2\pi}PV\int_{\R^2}\frac{-2(\dxi f(x,t)-\dxi f(x-y,t))}{[|y|^2+(f(x,t)-f(x-y,t))^2]^{3/2}}dy\\
&
-\frac{1}{2\pi}PV\int_{\R^2}\frac{
(\dxi f(x,t)-\dxi f(x-y,t))}{[|y|^2+(f(x,t)-f(x-y,t))^2]^{3/2}}B(x,y)dy,
\end{align*}
where
$$
B(x,y)=3\frac{|y|^2+(f(x,t)-f(x-y,t))\grad f(x-y)\cdot y}{|y|^2+(f(x)-f(x-y))^2}.
$$
Adding $I^i_2$ and $I^i_3$ one finds
\begin{equation}\label{i1i2}
I^i_2+I^i_3=-\frac{1}{2\pi}PV\int_{\R^2}\frac{\dxi f(x,t)-\dxi f(x-y,t)}{[|y|^2+(f(x,t)-f(x-y,t))^2]^{3/2}}C(x,y)dy,
\end{equation}
where $C=A+B-2.$

Consider
$$
M(t)=\max_{x\in\R^2}\left\{(\dxu f(x,t))^2+(\dxd f(x,t))^2 \right\}
=(\dxu f(x_t,t))^2+(\dxd f(x_t,t))^2.
$$ Next we follow the time derivative of $M(t)$ to find that $M'(t)\leq 0$ for almost every $t>0$ if $M(0)< 1/9$. This will yield the desired result.

 We obtain
$$M'(t)=2(\dxu f(x_t,t)\dxu f_t(x_t,t)+\dxd f(x_t,t)\dxd f_t(x_t,t))$$ for almost every $t$ (see \cite{DP2} for more details). It gives
$$
M'(t)=2\sum_{i=1,2}\dxi f (x_t,t)(I^i_2(x_t,t)+I^i_3(x_t,t)),
$$
due to the fact that at the maximum we have $$\dxu f(x_t,t) I^1_1(x_t,t)+\dxd f(x_t,t)I^2_1(x_t,t)=0.$$ Equation \eqref{i1i2} shows that it remains to check that $C(x_t,y)\geq 0.$ We write
$$
C(x,t)=1+3\frac{\Delta_y f(x)(\grad f(x)\cdot u- \Delta_y f(x))}{1+(\Delta_y f(x))^2},
$$
for $u=y/|y|$. It is easy to check that it is positive if $\|\grad f\|_{L^\infty}<1/3.$

\section{Global existence for initial data smaller than 1/3}\label{sec:init1}

Here we prove the existence of weak solutions for the Muskat contour equation. First we provide the notion of weak solution. It is possible to rewrite \eqref{ec} as follows:
\begin{equation}
f_t=\frac{\rho}{2\pi}\grad_x\cdot PV\int_{\R^2}\frac{y}{|y|^2}\frac{\Delta_yf(x)}{[1+(\Delta_yf(x))^2]^{1/2}}dy,
\label{muskatAT}
\end{equation}
where $\rho$ and $\Delta_yf(x)$ are defined as before. Then integrating by parts in the nonlinear term, it is easy to find that for any $\eta(x,t)\in C^\infty_c([0,T)\times\R^2)$, a weak solution $f$ should satisfy
\begin{multline}\label{weaksol}
\int_0^T\!\!\int_{\R^2}\eta_t(x,t)f(x,t)dxdt+\int_{\R^2}\eta(x,0)f_0(x)dx
\\
=\int_0^T\!\!\int_{\R^2}\grad_x\eta(x,t)\cdot\frac{\rho}{2\pi}PV\int_{\R^2}\frac{y}{|y|^2}
\frac{\Delta_yf(x)}{[1+(\Delta_yf(x))^2]^{1/2}}dy dxdt.
\end{multline}

 The main result we prove below is the following:
\begin{thm}\label{weakSOLthm}
Suppose that $\|f_0\|_{L^\infty}<\infty$ and $\|\grad f_0\|_{L^\infty}<1/3$.  Then there exists a weak solution of \eqref{weaksol} that satisfies
$$
f(x,t)\in C([0,T]\times\R^2)\cap L^\infty([0,T];W^{1,\infty}(\R^2)),
$$
for any $T>0$. In particular $f$ is a global in time Lipschitz continuous solution.
\end{thm}

We split the proof of Theorem \ref{weakSOLthm} in several sections.
A regularized model is defined below in \eqref{regularizedM} with solutions $f^\ep(x,t)$; here the model will be defined for a sufficiently small $\ep>0$. In Section \ref{sec:apriori} we prove some necessary a priori bounds for $f^\ep(x,t)$. They are used in Section \ref{sec:globalREG} to give global in time existence of classical solutions to the regularized model. Then, in Section \ref{secweakSOL} we explain how to obtain the weak solution as a limit as $\ep\to 0^+$; to this end we will establish to a strong convergence result.

 The regularized model is given by
\begin{align}
\begin{split}\label{regularizedM}
f^\ep_t(x,t)&=-\ep C\Lambda^{1-\ep}f^\ep+\ep \Delta f^\ep\\
&\quad+\frac{\rho}{2\pi}\grad_x\cdot PV\int_{\R^2} dy\frac{y}{|y|^{2-\ep}} ~ \frac{\Delta_yf^\ep(x)}{[1+(\Delta_yf^\ep(x))^2]^{1/2}},
\end{split}
\end{align}
where $C>0$ is an universal constant fixed below, the operator $\Lambda^{1-\ep}$ is a Fourier multiplier given by $\widehat{\Lambda^{1-\ep}f}(\xi)=|\xi|^{1-\ep}\widehat{f}(\xi)$ or equivalently using its integral from by
$$
\Lambda^{1-\ep}f(x)=c_\ep\int_{\R^2}\frac{f(x)-f(x-y)}{|y|^{3-\ep}}dy,
$$
with $\ep$ small enough. We define
$
\Delta f(x)=\dxu^2 f(x)+\dxd^2 f(x),
$
and $\Delta_yf(x)$ is given in \eqref{incre}.

In the next two subsections we write $f=f^\ep$ for the solution to \eqref{regularizedM} for the sake of simplicity of notation.

\subsection{A priori bounds}\label{sec:apriori}

For solutions of the regularized system \eqref{regularizedM} we get the following two a priori bounds $$\|f\|_{L^\infty}\leq \|f_0\|_{L^\infty},\qquad \|\grad f\|_{L^\infty}\leq \|\grad f_0\|_{L^\infty}<1/3.$$
The first one is obtained by checking the evolution of
$$
M(t)=\max_{x}f(x,t)=f(x_t,t).
$$
Here $x_t$ is thought of as the point where the maximum is attained.

For almost every $t$ we find
\begin{align*}
M'(t)=f_t(x_t,t)=-\ep C\Lambda^{1-\ep}f(x_t)+\ep \Delta f(x_t)+\frac{\rho}{2\pi}I(x_t),
\end{align*}
with
$$
I(x)=\grad_x\cdot PV\int_{\R^2}dy \frac{y}{|y|^{2-\ep}} ~ \frac{\Delta_yf(x)}{[1+(\Delta_yf(x))^2]^{1/2}}.
$$
Since
\begin{equation}
I(x)=\grad_x\cdot PV\int_{\R^2}dy \frac{x-y}{|x-y|^{2-\ep}} ~ \frac{\Delta_{x-y}f(x)}{[1+(\Delta_{x-y}f(x))^2]^{1/2}},
\label{Idefin}
\end{equation}
it is easy to find
\begin{align}
\begin{split}\label{nose}
I(x)&=\ep  PV\int_{\R^2}dy \frac{f(x)-f(y)}{|x-y|^{3-\ep}}  \frac1{[1+(\Delta_{x-y} f(x))^2]^{1/2}}\\
&\quad +PV\!\int_{\R^2}dy\frac{\grad f(x)\cdot(x-y)-(f(x)-f(y))}{|x-y|^{3-\ep}[1+(\Delta_{x-y} f(x))^2]^{3/2}}.
\end{split}
\end{align}
The previous formula shows that for $C$ large enough $$-\ep C\Lambda^{1-\ep} f(x_t)+\frac{\rho}{2\pi}I(x_t)\leq 0.$$ Then $M'(t)\leq 0$ for a.e. $t\in(0,T]$ because $\Delta f(x_t)\leq 0$ and therefore $M(t)\leq M(0)$. Analogously $m(t)\geq m(0)$.

Next we consider the evolution of
$$
L(t)=\max_{x\in \R^2}(\dxu f(x,t))^2+(\dxd f(x,t))^2=(\dxu f(x'_t,t))^2+(\dxd f(x'_t,t))^2.
$$
We can proceed as in the previous section, but in this case more terms will appear. In $\dxi f_t$ we have analogous terms that can be handled as before. Terms with the correct sign, that appear due to $-\ep\Lambda^{1-\ep}f$ and $\ep\Delta f$ in \eqref{regularizedM}.  And a new element $J^i(x)$ has terms which are given by
$$
J^i(x)=\ep \int_{\R^2}dy\frac{\dxi f(x)-\dxi f(x-y)}{|y|^{3-\ep}}\frac{1}{[1+(\Delta_y f(x))^2]^{3/2}}.
$$
That is
$$
\dxi f_t(x)=-\ep C\Lambda^{1-\ep}\dxi f(x)+\ep\Delta\dxi f(x)+\frac{\rho}{2\pi}J^i(x)+\mbox{``Analogous terms"}.
$$
In checking the time derivative of $L$
$$L'(t)=2\sum_{i=1,2}\dxi f(x'_t,t)\dxi f_t(x'_t,t),$$
all the terms are handled as before but for
$$
\frac{\rho}{2\pi}\sum_{i=1,2}\dxi f(x'_t)J^i(x'_t).
$$
But at this point it is easy to check that
$$
-\ep C\sum_{i=1,2}\dxi f(x'_t)\Lambda^{1-\ep}\dxi f(x'_t)+ \frac{\rho}{2\pi}\sum_{i=1,2}\dxi f(x'_t)J^i(x'_t)\leq 0,
$$
for $C$ big enough. Therefore $L'(t)\leq 0$ if $\sqrt{L(t)}<1/3$ for almost every $t$. This yields the desired maximum principle.

\subsection{Global existence for the regularized model}\label{sec:globalREG}

We consider regular initial data $f_0\in H^4$ for the system \eqref{regularizedM}. Local existence can easily be proved using the energy method following the arguments for the non-regularized Muskat problem \eqref{ec}, as in \cite{DP}.

As we did for \eqref{ec}, it follows that
\begin{align*}
\frac{d}{dt}\|f\|_{L^2}(t)&=-\frac{\rho}{\pi}\int_{\R^2}\int_{\R^2}
\frac{1+\ep}{|x-y|^{1-\ep}}(1-[1+(\Delta_y f(x))^2]^{-1/2}) dx dy\\
&\quad -2C\ep\|\Lambda^{(1-\ep)/2} f\|_{L^2}(t)-2\ep\|\grad f\|_{L^2}(t).
\end{align*}
Therefore
$
\|f\|_{L^2}(t)\leq \|f_0\|_{L^2}.
$

\begin{rem}
The global existence theorem for weak solutions can also be found with
$$
\|f_0\|_{L^2}<\infty\quad \mbox{instead of} \quad \|f_0\|_{L^\infty}<\infty.
$$
We chose the version above because it is more general.  We see that if the solution satisfies initially a $L^2$ bound then
$
f(x,t)\in L^\infty([0,T];L^2(\R^2)).
$
\end{rem}

Next, we consider the evolution of
\begin{align*}
\int_{\R^2} \dxu^3f\dxu^3f_t dx&\leq -C\ep\|\Lambda^{(1-\ep)/2} \dxu^3f\|_{L^2}(t)
-\ep\|\grad\dxu^3f\|^2_{L^2}+I_1+I_2,
\end{align*}
where
$$
I_1=\frac{\rho}{2\pi}\int_{\R^2} \dxu^3f(x)\dxu^3\Big(PV\int_{\R^2}
\frac{(\grad f(x)-\grad f(x-y))\cdot y}{|y|^{3-\ep}}dy\Big) dx,
$$
\begin{align*}
I_2=\frac{\rho}{2\pi}\int_{\R^2}\dxu^3f(x)\dxu^3\Big(PV\int_{\R^2} &\frac{(\grad f(x)-\grad f(x-y))\cdot y}{|y|^{3-\ep}}\frac{1-[1+(\Delta_y f(x))^2]^{3/2}}{[1+(\Delta_y f(x))^2]^{3/2}}dy \Big) dx.
\end{align*}
The term $\grad f(x)$ cancels out in $I_1$ due to the PV and an integration by parts shows that
$$
I_1=-\frac{\rho}{2\pi}C(\ep)\int_{\R^2} \dxu^3f(x)\Lambda^{1-\ep}\dxu^3f(x)dx\leq 0.$$
For $I_2$ one finds
\begin{align*}
I_2=\frac{\rho}{2\pi}\int_{\R^2}\dxu^4f(x)\dxu^2\Big(PV\int_{\R^2} \frac{(\grad f(x)-\grad f(x-y))\cdot y}{|y|^{3-\ep}}\frac{[1+(\Delta_y f(x))^2]^{3/2}-1}{[1+(\Delta_y f(x))^2]^{3/2}}dy \Big) dx,
\end{align*}and the splitting $I_2=J_1+J_2+J_3+J_4$ gives
\begin{align*}
J_1=\frac{\rho}{2\pi} \int_{\R^2}\dxu^4f(x)\int_{\R^2} \frac{(\grad\dxu^2f(x)-\grad\dxu^2f(x-y))\cdot y}{|y|^{3-\ep}}\frac{[1+(\Delta_y f(x))^2]^{3/2}-1}{[1+(\Delta_y f(x))^2]^{3/2}}dy  dx,
\end{align*}
\begin{align*}
J_2=\frac{3\rho}{\pi} \int_{\R^2}\dxu^4f(x)\int_{\R^2}
\frac{(\grad\dxu f(x)-\grad\dxu f(x-y))\cdot y}{|y|^{3-\ep}}\frac{\Delta_y\dxu f(x)\Delta_y f(x)}{[1+(\Delta_y f(x))^2]^{5/2}}dy  dx,
\end{align*}
\begin{align*}
J_3=\frac{3\rho}{2\pi} \int_{\R^2}\dxu^4f(x)\int_{\R^2} \frac{(\grad f(x)-\grad f(x-y))\cdot y}{|y|^{3-\ep}}\frac{\Delta_y\dxu^2 f(x)\Delta_y f(x)}{[1+(\Delta_y f(x))^2]^{5/2}}dy  dx,
\end{align*}
\begin{align*}
J_4=\frac{3\rho}{2\pi} \int_{\R^2}\dxu^4f(x)\int_{\R^2} \frac{(\grad f(x)-\grad f(x-y))\cdot y}{|y|^{3-\ep}}
(\Delta_y\dxu f(x))^2\frac{1-4(\Delta_y f(x))^2}{[1+(\Delta_y f(x))^2]^{7/2}}dy  dx.
\end{align*}

For $J_1$ we proceed as follows
\begin{multline*}
\left| J_1 \right| =\frac{\rho}{2\pi} \left(\int_{|y|>1} dx\int_{\R^2} dx+\int_{|y|<1} dy\int_{\R^2} dx\right)
\leq C(\ep)(\|f\|_{L^\infty}+1)\|\grad\dxu^2f\|_{L^2}\|\dxu^4f\|_{L^2}.
\end{multline*}
The identity
$$
\dxi\dxu f(x)-\dxi\dxu f(x-y)=\int_0^1\grad\dxi\dxu f(x+(s-1)y)\cdot y  ds,
$$
yields
\begin{multline*}
\left| J_2 \right|
\leq
\frac{3\rho}{\pi} \int_0^1ds\int_{|y|<1} \frac{dy}{|y|^{2-\ep}}\int_{\R^2} dx
|\partial_x^4f(x)||\grad^2\dxu f(x+(s-1)y)|
(|\dxu f(x)|+|\dxu f(x-y)|)
\\
+
\frac{3\rho}{\pi} \int_0^1\!\!ds\int_{|y|>1}~ \frac{dy}{|y|^{3-\ep}}\int_{\R^2} \!\!dx ~ |\dxu^4f(x)||\grad^2\dxu f(x\!+\!(s\!-\!1)y)|
\\
\times(|\dxu f(x)|\!+\!|\dxu f(x\!-\!y)|)(|f(x)|\!+\!|f(x\!-\!y)|),
\end{multline*}
and therefore
$$
\left| J_2 \right| \leq C(\ep)(1+\|f\|_{L^\infty})\|\grad f\|_{L^\infty}\|\dxu^4f\|_{L^2}\|\grad^2\dxu f\|_{L^2}.
$$

In $J_3$ we use the splitting $J_3=K_1+K_2$ where
$$
K_1=\frac{3\rho}{2\pi}\int_{|y|>1} dy\int_{\R^2} dx,\qquad
K_2=\frac{3\rho}{2\pi}\int_{|y|<1} dy\int_{\R^2} dx,
$$
and then
\begin{align*}
\left| K_1 \right|
&\leq \frac{3\rho}{\pi}\|\grad f\|_{L^\infty} \int_{|y|>1} \frac{dy}{|y|^{3-\ep}}\int_{\R^2} dx
|\dxu^4f(x)|(|\dxu^2f(x)|+|\dxu^2f(x-y)|) \\
&\leq C\|\grad f\|_{L^\infty}\|\dxu^4f\|_{L^2}\|\dxu^2f\|_{L^2}\\
&\leq C\|\grad f\|_{L^\infty}\|\dxu^4f\|_{L^2}(\|f\|_{L^2}+\|\dxu^3f\|_{L^2}).
\end{align*}
The equality
$$
\dxu^2f(x)-\dxu^2f(x-y)=\int_0^1\grad\dxu^2f(x+(s-1)y)\cdot y ds,
$$
allows us to obtain
\begin{align*}
\left| K_2 \right|
&\leq \frac{3\rho}{\pi}\|\grad f\|_{L^\infty}\!\!\int_0^1\!\!ds \int_{|y|<1}~ \frac{dy}{|y|^{2-\ep}}\int_{\R^2} dx|\partial_x^4f(x)|
|\grad\dxu^2f(x+(s-1)y)|\\
&\leq C\|\grad f\|_{L^\infty}\|\dxu^4f\|_{L^2}\|\grad\dxu^2f\|_{L^2}.
\end{align*}
In $J_4$ we use the splitting $J_4=K_3+K_4$ where
$$
K_3=\frac{3\rho}{2\pi}\int_{|y|>1} dy\int_{\R^2} dx,\qquad
K_4=\frac{3\rho}{2\pi}\int_{|y|<1} dy\int_{\R^2} dx,
$$
and then
\begin{align*}
\left| K_3 \right|
&\leq C\|\grad f\|^2_{L^\infty} \int_{|y|>1} \frac{dy}{|y|^{4-\ep}}\int_{\R^2} dx
|\dxu^4f(x)|(|\dxu f(x)|+|\dxu f(x-y)|) \\
&\leq C\|\grad f\|^2_{L^\infty}\|\dxu^4f\|_{L^2}\|\dxu f\|_{L^2}\\
&\leq C\|\grad f\|^2_{L^\infty}\|\dxu^4f\|_{L^2}(\|f\|_{L^2}+\|\dxu^3f\|_{L^2}).
\end{align*}
The equality
$$
\dxu^2f(x)-\dxu^2f(x-y)=\int_0^1\grad\dxu^2f(x+(s-1)y)\cdot y ds,
$$
allows us to obtain
\begin{align*}
\left| K_4 \right|
&\leq C\|\grad f\|_{L^\infty}\int_0^1\!ds\int_0^1\!dr \int_{|y|<1} \!|y|^{\ep-2}dy\int_{\R^2}\! dx |\dxu^4f(x)|\\
&\qquad\qquad\qquad\qquad\times (|\grad\dxu f(x+(s-1)y)||\grad\dxu f(x+(r-1)y)|) \\
&\leq C\|\grad f\|_{L^\infty}\|\dxu^4f\|_{L^2}\|\grad\dxu f\|^2_{L^4}.
\end{align*}
The following estimate
\begin{equation*}
\|\dxi\dxu f\|_{L^4}^4\leq 3\|\grad f\|_{L^\infty}\|\dxi\dxu f\|_{L^4}^2\|\dxi\dxu^2f\|_{L^2},
\end{equation*}
yields
$$
\left| K_4 \right| \leq C\|\grad f\|^2_{L^\infty}\|\dxu^4f\|_{L^2}\|\grad\dxu^2 f\|_{L^2}.
$$
Using Young's inequality
\begin{align*}
\frac{d}{dt}\|\dxu^3f\|^2_{L^2}\leq C(\ep)(\|f\|^2_{L^\infty}+\|\grad f\|^2_{L^\infty}&+\|\grad f\|^4_{L^\infty}+1)\|f\|^2_{H^3}.
\end{align*}Proceeding in a similar manner, at this point it is easy to find
\begin{align*}
\frac{d}{dt}\|\dxd^3f\|^2_{L^2}\leq C(\ep)(\|f\|^2_{L^\infty}+\|\grad f\|^2_{L^\infty}&+\|\grad f\|^4_{L^\infty}+1)\|f\|^2_{H^3},
\end{align*}
and therefore
\begin{align*}
\frac{d}{dt}\|f\|^2_{H^3}\leq C(\ep)(\|f\|^2_{L^\infty}+\|\grad f\|^2_{L^\infty}&+\|\grad f\|^4_{L^\infty}+1)\|f\|^2_{H^3}.
\end{align*}
The Gronwall inequality then yields
$$
\|f\|^2_{H^3}(t)\leq \|f_0\|^2_{H^3}\exp\Big(\int_0^t C(\ep)G(s)ds\Big),
$$
for
$$
G(s)=(\|f\|^2_{L^\infty}(s)+\|\grad f\|^2_{L^\infty}(s)+\|\grad f\|^4_{L^\infty}(s)+1).
$$ We find $f\in C([0,T];H^3(\R))$ for any $T>0$ by the a priori bounds.

  For the argument in next sections we will need $f\in C([0,T];H^4(\R))$ for any $T>0$.  Therefore we consider the evolution of four derivatives. Most of the terms can be controlled as before. We will show how to deal with the rest using the estimate of the $H^3$ norm. Since
\begin{align*}
\int_{\R^2} \dxu^4f\dxu^4f_t dx&\leq -C\ep\|\Lambda^{(1-\ep)/2} \dxu^4f\|_{L^2}(t)
-\ep\|\grad\dxu^4f\|^2_{L^2}+L_1+L_2,
\end{align*}
where
$$
L_1=\frac{\rho}{2\pi}\int_{\R^2} \dxu^4f(x)\dxu^4\Big(PV\int_{\R^2}
\frac{(\grad f(x)-\grad f(x-y))\cdot y}{|y|^{3-\ep}}dy\Big) dx,
$$
\begin{align*}
L_2=\frac{\rho}{2\pi}\int_{\R^2}\dxu^4f(x)\dxu^4\Big(PV\int_{\R^2} &\frac{(\grad f(x)-\grad f(x-y))\cdot y}{|y|^{3-\ep}}\frac{1-[1+(\Delta_y f(x))^2]^{3/2}}{[1+(\Delta_y f(x))^2]^{3/2}}dy \Big) dx.
\end{align*}
The term $L_1$ has the correct sign as $I_1$. For $L_2$ one finds
\begin{align*}
L_2=\frac{\rho}{2\pi}\int_{\R^2}\dxu^5f(x)\dxu^3\Big(PV\int_{\R^2} &\frac{(\grad f(x)-\grad f(x-y))\cdot y}{|y|^{3-\ep}}
\frac{[1+(\Delta_y f(x))^2]^{3/2}-1}{[1+(\Delta_y f(x))^2]^{3/2}}dy \Big) dx,
\end{align*}and the splitting $L_2=M_1+M_2+M_3+M_4$ gives
\begin{align*}
M_1=\frac{\rho}{2\pi} \int_{\R^2}\dxu^5f(x)\int_{\R^2}\frac{(\grad\dxu^3f(x)-\grad\dxu^3f(x-y))\cdot y}{|y|^{3-\ep}}\frac{[1+(\Delta_y f(x))^2]^{3/2}-1}{[1+(\Delta_y f(x))^2]^{3/2}}dy  dx,
\end{align*}
\begin{align*}
M_2=\frac{9\rho}{2\pi} \int_{\R^2}\dxu^5f(x)\int_{\R^2} \frac{(\grad\dxu^2 f(x)-\grad\dxu^2 f(x-y))\cdot y}{|y|^{3-\ep}}\frac{\Delta_y\dxu f(x)\Delta_y f(x)}{[1+(\Delta_y f(x))^2]^{5/2}}dy  dx,
\end{align*}
\begin{align*}
M_3=\frac{9\rho}{2\pi} \int_{\R^2}\dxu^5f(x)\int_{\R^2}
\frac{(\grad\dxu f(x)-\grad\dxu f(x-y))\cdot y}{|y|^{3-\ep}}\dxu\Big(\frac{\Delta_y\dxu f(x)\Delta_y f(x)}{[1+(\Delta_y f(x))^2]^{5/2}}\Big)dy  dx,
\end{align*}
\begin{align*}
M_4=\frac{3\rho}{2\pi} \int_{\R^2}\dxu^5f(x)\int_{\R^2}
\frac{(\grad f(x)-\grad f(x-y))\cdot y}{|y|^{3-\ep}}\dxu^2\Big(\frac{\Delta_y\dxu f(x)\Delta_y f(x)}{[1+(\Delta_y f(x))^2]^{5/2}}\Big)dy  dx.
\end{align*}
For $M_1$ and $M_2$ we obtain as before
$$
|M_1|+|M_2|\leq C(\ep)(1+\|f\|_{L^\infty})(\|\grad f\|_{L^\infty}+1)\|\dxu^5f\|_{L^2}\|f\|_{H^4}.
$$
In $M_3$ we use the splitting $M_3=N_1+N_2$ where
\begin{align*}
N_1=\frac{9\rho}{2\pi} \int_{\R^2}\dxu^5f(x)\int_{\R^2}
\frac{(\grad\dxu f(x)-\grad\dxu f(x-y))\cdot y}{|y|^{3-\ep}}\Delta_y\dxu^2 f(x)\frac{\Delta_y f(x)}{[1+(\Delta_y f(x))^2]^{5/2}}dy  dx,
\end{align*}
\begin{align*}
N_2=\frac{9\rho}{2\pi} \int_{\R^2}\dxu^5f(x)\int_{\R^2} &
\frac{(\grad\dxu f(x)-\grad\dxu f(x-y))\cdot y}{|y|^{3-\ep}}\\
&\times(\Delta_y\dxu f(x))^2\frac{1-4(\Delta_y f(x))^2}{[1+(\Delta_y f(x))^2]^{7/2}}\Big)dy  dx,
\end{align*}
We take
$$N_1=\frac{9\rho}{2\pi}\int_{|y|>1} dy\int_{\R^2} dx+\frac{9\rho}{2\pi}\int_{|y|<1} dy\int_{\R^2} dx,
$$
to find as before
$$
|N_1|\leq C(\ep)\|\dxu^5f\|_{L^2}\|\grad\dxu f\|_{L^4}(\|\dxu^2f\|_{L^4}+\|\grad\dxu^2 f\|_{L^4}).
$$
By Sobolev embedding
$$
|N_1|\leq C(\ep)\|\dxu^5f\|_{L^2}\|f\|_{H^3}(\|f\|_{H^3}+\|f \|_{H^4}).
$$
Similarly for $N_2$
\begin{align*}N_2=\frac{9\rho}{2\pi}\int_{|y|>1} dy\int_{\R^2} dx+\frac{9\rho}{2\pi}\int_{|y|<1} dy\int_{\R^2} dx,
\end{align*}
therefore
\begin{align*}
|N_2|\leq C(\ep)\|\dxu^5f\|_{L^2}\|\dxu f\|_{L^\infty}(& \|\dxu f\|_{L^\infty}\|\grad\dxu f\|_{L^2}\\
&+\|\grad^2\dxu f\|_{L^4}\|\grad\dxu f\|_{L^4}),
\end{align*}
which yields
\begin{align*}
|N_2|\leq C(\ep)\|\dxu^5f\|_{L^2}\|\grad f
_0\|_{L^\infty}\|f\|_{H^3}(\|\grad f_0\|_{L^\infty}+\| f\|_{H^4}),
\end{align*}
For $M_4$ we split further $M_4=N_3+N_4+N_5$
\begin{align*}
N_3=\frac{3\rho}{2\pi} \int_{\R^2}\dxu^5f(x)\int_{\R^2}
\frac{(\grad f(x)-\grad f(x-y))\cdot y}{|y|^{3-\ep}}\Delta_y\dxu^3 f(x)\frac{\Delta_y f(x)}{[1+(\Delta_y f(x))^2]^{5/2}}dy  dx,
\end{align*}
\begin{align*}
N_4=\frac{9\rho}{2\pi} \int_{\R^2}\dxu^5f(x)\int_{\R^2} &
\frac{(\grad f(x)-\grad f(x-y))\cdot y}{|y|^{3-\ep}}\\
&\times \Delta_y\dxu^2 f(x)\Delta_y\dxu f(x)\frac{1-4
(\Delta_y f(x))^2}{[1+(\Delta_y f(x))^2]^{7/2}}dy  dx,
\end{align*}
\begin{align*}
N_5=\frac{3\rho}{2\pi} \int_{\R^2}\dxu^5f(x)\int_{\R^2} &
\frac{(\grad f(x)-\grad f(x-y))\cdot y}{|y|^{3-\ep}}\\
&\times (\Delta_y\dxu f(x))^3\frac{20(\Delta_y f(x))^3\!-\!8\Delta_y f(x)\!-\!7}{[1+(\Delta_y f(x))^2]^{9/2}}dy  dx,
\end{align*}
For $N_3$ one finds
\begin{align*}
|N_3|&\leq C(\ep)\|\dxu^5f\|_{L^2}\|\grad f\|_{L^\infty}
(\|\dxu^3 f\|_{L^2}+\|\grad\dxu^3 f\|_{L^2})\\
&\leq C(\ep)\|\dxu^5f\|_{L^2}\|\grad f_0\|_{L^\infty}(\|f\|_{H^3}+\|f\|_{H^4}),
\end{align*}
and similarly  for $N_4$
\begin{align*}
|N_4|&\leq C(\ep)\|\dxu^5f\|_{L^2}\|\grad f\|_{L^\infty}
(\|\dxu^2 f\|_{L^2}\|\dxu f\|_{L^\infty}+\|\grad\dxu^2 f\|_{L^4}\|\grad\dxu f\|_{L^4})\\
&\leq C(\ep)\|\dxu^5f\|_{L^2}\|\grad f_0\|_{L^\infty}\|f\|_{H^3}(\|\grad f_0\|_{L^\infty}+\|f\|_{H^4}).
\end{align*}
Finally, for $N_5$ we conclude that
\begin{align*}
|N_5|&\leq C(\ep)\|\dxu^5f\|_{L^2}\|\grad f\|_{L^\infty}
(\|\dxu f\|_{L^2}\|\dxu f\|^2_{L^\infty}+\|\grad\dxu f\|^3_{L^6})\\
&\leq C(\ep)\|\dxu^5f\|_{L^2}\|\grad f_0\|_{L^\infty}(\|f\|_{H^3}\|\grad f_0\|^2_{L^\infty}+\|f\|^3_{H^3}),
\end{align*}
by Sobolev embedding.

If we gather all the estimates above and use Young's inequality, it is not difficult to check that
\begin{align*}
\int_{\R^2}\dxu^4f(x)\dxu^4f_t(x)dx\leq& C(\ep)(1\!+\!\|f_0\|^2_{L^\infty})(1\!+\!\|\grad f_0\|^2_{L^\infty})(1\!+\!\|f\|^2_{H^3})\|f\|^2_{H^4}\\
&+C(\ep)(1+\|\grad f_0\|^6_{L^\infty})(1+\|f\|^6_{H^3}).
\end{align*}
A repetition of the argument for $\dxd^4$ gives
\begin{align*}
\int_{\R^2}\dxd^4f(x)\dxu^4f_t(x)dx\leq& C(\ep)(1\!+\!\|f_0\|^2_{L^\infty})(1\!+\!\|\grad f_0\|^2_{L^\infty})(1\!+\!\|f\|^2_{H^3})\|f\|^2_{H^4}\\
&+C(\ep)(1+\|\grad f_0\|^6_{L^\infty})(1+\|f\|^6_{H^3}).
\end{align*}
Therefore
\begin{align*}
\frac{d}{dt}\|f\|^2_{H^4}\leq& C(\ep)(1\!+\!\|f_0\|^2_{L^\infty})(1\!+\!\|\grad f_0\|^2_{L^\infty})(1\!+\!\|f\|^2_{H^3})\|f\|^2_{H^4}\\
&+C(\ep)(1+\|\grad f_0\|^6_{L^\infty})(1+\|f\|^6_{H^3}).
\end{align*}
We use the Gronwall inequality and additionally the control of the $H^3$ norm to obtain the desired global estimate for $H^4$.

\subsection{Taking $\ep\to 0^+$}\label{secweakSOL}

This section ends the proof of Theorem \ref{weakSOLthm} by showing that solutions of the regularized system converge to a weak solution.

First we approximate the initial data to have a global solution of the regularized system. An approximation to the identity $\zeta \in C^\infty_c(\R^2)$ is defined as follows:
\begin{equation}\label{zeta}
\int_{\R^2} dx ~ \zeta(x) = 1,
\quad
\zeta \ge 0,\quad \zeta(x)=\zeta(-x),
\quad\mbox{where}\quad \zeta_\ep(x) = \zeta(x/\ep) /\ep^2.
\end{equation}
Then, for any $f_0\in W^{1,\infty}(\R^2)$ and $\|\grad f_0\|_{L^\infty}<1/3$, we define the initial data for the regularized system as follows
$$
f_0^\ep (x) = \frac{(\zeta_\ep * f_0)(x)}{1+\ep^2 |x|^2}.
$$
Notice that $f_0^\ep \in H^s(\R)$ for any $s>0$, and
$
\|f_0^\ep\|_{L^\infty}\leq \|f_0\|_{L^\infty}.
$
More importantly,
$
\|\grad f_0^\ep\|_{L^\infty}<1/3
$
if $\ep$ is sufficiently small ($\ep$ depends upon the size of $\|f_0\|_{L^\infty}$).
Therefore global existence of the regularized system \eqref{regularizedM} holds with initial data $f_0^\ep$ under the condition that $\ep>0$ is small enough.

Now consider the solutions $\{f^\ep\}$ to the regularized system
\eqref{regularizedM} with initial data given by the $f_0^\ep$ as described above.
Integration by parts provides
\begin{align}
\begin{split}\label{weaksolep}
&\int_0^T\!\!\int_{\R^2}\eta_tf^\ep dxdt+\int_{\R^2}\ep(C\Lambda^{1-\ep}-\Delta)\eta\, f^\ep dxdt+\int_{\R^2}\eta(x,0)f^\ep_0(x)dx\\
&=\int_0^T\!\!\int_{\R^2} \grad_x\eta(x,t)\cdot\frac{\rho}{2\pi}PV\int_{\R^2}\frac{y}{|y|^{2-\ep}}
\frac{\Delta_yf^\ep(x)}{[1+(\Delta_yf^\ep(x))^2]^{1/2}}dy dxdt,
\end{split}
\end{align}
for any $\eta\in C_c^\infty([0,T)\times\R^2)$.

Now we send $\ep\to0^+$ to in order to obtain \eqref{weaksol}. The third integral above converges as a result of the properties of the the approximation to the identity which was previously introduced. The second integral converges to $0$ because of the bound $\|f^\ep\|_{L^\infty}(t)\leq\|f_0\|_{L^\infty}$. Together with the other bound ($\|\grad f^\ep\|_{L^\infty}(t) < 1/3$), we find the existence of a subsequence (denoted again by $f^\ep$) that converges in the weak* topology to a function $f\in L^\infty([0,T];W^{1,\infty}(\R^2))$ by the Banach-Alaoglu theorem. This provides the solution $f$ and implies the convergence of the first integral in \eqref{weaksolep}. It remains to check that as $\ep \to 0^+$ we have
\begin{multline*}
\int_0^T ~ dt ~
\int_{\R^2} ~ dx ~
\grad\eta(x,t)\cdot\frac{\rho}{2\pi}PV\!\!\int_{\R^2}~ dy~ \frac{y}{|y|^{2-\ep}}\frac{\Delta_yf^\ep(x)}{[1+(\Delta_yf^\ep(x))^2]^{1/2}}
\\
\to
\int_0^T ~ dt ~
\int_{\R} ~ dx ~
\grad\eta(x,t)\cdot\frac{\rho}{2\pi}PV\!\!\int_{\R^2}~ dy~ \frac{y}{|y|^2}\frac{\Delta_yf(x)}{[1+(\Delta_yf(x))^2]^{1/2}}.
\end{multline*}

We let $B_R$ denote the open ball of radius $R$ and center $(0,0)$, then we {\it claim} that there is a subsequence (denoted again by $f^\ep$) such that
\begin{equation}\label{claim.conv}
\|f^\ep-f\|_{L^\infty([0,T]\times B_R)}\to 0,\quad \mbox{as}\quad\ep\to 0.
\end{equation}
We will prove this at the end of the section by using a strong convergence theorem. Since $f^\ep \in C([0,T]\times \R)$ for any $\ep>0$ and, up to a subsequence, $f^\ep$ converges to $f$ on compact sets, we obtain $f\in C([0,T]\times \R)$.

Choose $M>0$ so that  supp$(\eta)\subseteq B_M$. For any small $\delta >0$ and any large $L \gg 1$, with $L>M+1$ we split the integral as
\begin{equation}
\label{split.integral}
\int_{\R^2}~ dy~
=
\int_{B_\delta}~ dy~
+
\int_{B_L - B_\delta}~ dy~
+
\int_{B_L^c}~ dy.
\end{equation}
The first and last integrals separately are arbitrarily small independent of $\ep$ for $L>0$ sufficiently large and for $\delta>0$ sufficiently small:
The bound
$$
\left|
\frac{\Delta_yf^\ep(x)}{[1+(\Delta_yf^\ep(x))^2]^{1/2}}  \right|\le 1.
$$
yields
\begin{multline*}
\left|
\int_0^T dt
\int_{\R^2} dx ~\grad\eta(x,t)\cdot\frac{\rho}{2\pi}PV\int_{B_\delta}~ dy~ \frac{y}{|y|^{2-\ep}}\frac{\Delta_yf^\ep(x)}{[1+(\Delta_yf^\ep(x))^2]^{1/2}}\right|
\le \rho \|\grad \eta\|_{L^1([0,T]\times \R^2)}\delta.
\end{multline*}
 For the integral on $B_L^c$ we note that
$$
\frac{z}{[1+z^2]^{1/2}}= \int_0^1\frac{d}{ds} \frac{sz}{[1+(sz)^2]^{1/2}}ds=z \int_0^1\frac{1}{[1+(sz)^2]^{3/2}} ds,
$$
and therefore
$$
\frac{z}{[1+z^2]^{1/2}}= z\Big(1+ \int_0^1\frac{1-[1+(sz)^2]^{3/2}}{[1+(sz)^2]^{3/2}} ds\Big)
=
z\left(1- z^2\int_0^1s^2 h(sz) ds\right),
$$
where $h(sz)=(3+3(sz)^2+(sz)^4)/([1+(sz)^2]^{3/2}(1+[1+(sz)^2]^{3/2}))$.
This expression allows us to split
\begin{multline*}
PV\!\!\int_{B_L^c}~ dy~\frac{y_i}{|y|^{2-\ep}} \frac{\Delta_yf^\ep(x)}{[1+(\Delta_yf^\ep(x))^2]^{1/2}}=
-R_i^{\ep,L}(f^\ep)
\\
-PV\int_{B_L^c}~ dy~\frac{y_i}{|y|^{2-\ep}} (\Delta_yf^\ep(x))^3 \int_0^1s^2h(s\Delta_y f^\ep(x))ds,
\end{multline*}
for $i=1,2$. Here $R_i^{\ep,L}$ has the form
$$
R_i^{\ep,L}(f^\ep) \eqdef PV\!\!\int_{B_L^c} dy~\frac{f^\ep(x-y) y_i}{|y|^{3-\ep}},
$$
with the principal value at infinity.
   On the other hand in the second term in the left hand side the principal value is not necessary and we obtain
$$
\left|
\int_{B_L^c} dy \int_0^1 ds
\right|
\le
C\| f^\ep \|^3_{L^\infty}
\int_{B_L^c}~ \frac{dy}{|y|^{4-\ep}}\le
\frac{C\| f_0 \|^3_\infty}{L}.
$$
It remains to show a similar bound for
$$
I^L\eqdef \int_{B_M} dx ~\grad \eta(x,t)\cdot(R_1^{\ep,L}(f^\ep),R_2^{\ep,L}(f^\ep)).
$$
The principal value yields $I^L=\D\lim_{n\to\infty}I^L_n$ where
$$
I^L_n=\int_{B_M} dx~ \grad \eta(x,t)\cdot\int_{B_n\setminus B_L}dy~ f^\ep(x-y)\frac{y}{|y|^{3-\ep}}.
$$
Integration by parts provides
\begin{align*}
I_n^L=\int_{B_M}dx~\eta(x,t)\Big(&\int_0^{2\pi}\Big(\frac{f(x-nu)}{n^{1-\ep}}-\frac{f(x-Lu)}{L^{1-\ep}}\Big) d\theta\\
&+(1-\ep)\int_{B_n\setminus B_L}\frac{f^\ep(x-y)}{|y|^{3-\ep}}dy\Big),
\end{align*}
for $|u|=1$, which allows us to bound the following term as
$$
|I_n^L|\leq C\|\eta\|_{L^1}\|f\|_{L^\infty}\left(\frac{1}{L^{1-\ep}}+\frac{1}{n^{1-\ep}}\right).
$$
Hence
$$
|I^L|\leq C\|\eta\|_{L^1}\|f_0\|_{L^\infty}/L^{1/2}
$$
and we conclude that $I^L$ is arbitrarily small if $L$ is arbitrarily large.

For the last integral we recall that we have uniform convergence on compact sets. Due to $y \in B_L - B_\delta$ and $x\in B_M$ we have
\begin{multline*}
\int_0^T ~ dt ~
\int_{\R^2} ~ dx ~
\grad\eta(x,t)\cdot\frac{\rho}{2\pi}PV\!\!\int_{B_L\setminus B_\delta}
~ dy~ \frac{y}{|y|^{2-\ep}}\frac{\Delta_yf^\ep(x)}{[1+(\Delta_yf^\ep(x))^2]^{1/2}}
\\
\to
\int_0^T ~ dt ~
\int_{\R} ~ dx ~
\grad\eta(x,t)\cdot\frac{\rho}{2\pi}PV\!\!\int_{B_L\setminus B_\delta}~ dy~ \frac{y}{|y|^2}\frac{\Delta_yf(x)}{[1+(\Delta_yf(x))^2]^{1/2}},
\end{multline*}
as $\ep\to0^+$.

For $L$ sufficiently large and $\delta>0$ sufficiently small, we conclude by taking $\ep \to 0^+$.

It remains to prove the strong convergence in $L^\infty([0,T];L^\infty(B_R))$ for any $R>0$ which was claimed in \eqref{claim.conv}. The idea is to use the weak space $W_*^{-2,\infty}(B_R)$ to obtain bounds for $f_t^\ep(x,t)$ which are uniform:
\begin{equation}
\begin{split}
\sup_{t\in [0,T]}\| f_t^\ep\|_{W_*^{-2,\infty}(B_{R})}(t)
\le C  \| f_0 \|_{L^\infty(\R^2)},
\end{split}
\label{timeBOUND}
\end{equation}
where $C$ does not depend on $R$ or $\ep$. For $v\in L^\infty(B_{R})$ we consider the norm $\|\cdot\|_{W_*^{-2,\infty}(B_{R})}$ as follows:
$$
\| v\|_{W_*^{-2,\infty}(B_{R})}=  \sup_{\phi \in W^{2,1}_0(B_{R})\,:\, \| \phi \|_{2,1} \le 1} \left| \int_{B_R}\phi(x) v(x) dx \right|,
$$
where $W^{2,1}_0(B_{R})=\overline{C_c^\infty(B_{R})}^{W^{2,1}}$. Now the Banach space $W_*^{-2,\infty}(B_{R})$ is defined to be the completion of $L^\infty(B_{R})$ with respect to this norm $\| \cdot \|_{W_*^{-2,\infty}(B_{R})}$. We have the following result for convergence in this space (see \cite{PDPB} Lemma 4.3):

\begin{lemma}\label{convergence}
Consider a sequence $\{u_m\}$ in $C([0,T]\times B_R)$ that is uniformly bounded in the space
$L^\infty ([0,T]; W^{1,\infty}(B_R))$.
Assume further that the weak derivative $\partial_t u_m$ is in $L^\infty([0,T];L^\infty(B_R))$ (not necessarily uniform) and is uniformly bounded in $L^\infty ([0,T];W_*^{-2,\infty}(B_{R}))$. Finally suppose that $\partial_{x_i} u_m\in C([0,T]\times B_R)$ for $i=1,$ $2$ and any $m$ (not necessarily uniform). Then there exists a subsequence of $u_m$ that converges strongly in
$L^\infty ([0,T];L^\infty(B_R))$.
\end{lemma}

By applying this lemma the strong convergence {\it claimed} in \eqref{claim.conv} is obtained. It only remains to check the hypothesis of the lemma. For any regularized solution $f^\ep$ to \eqref{regularizedM} we need $f_t^\ep$ in $L^\infty([0,T];L^\infty(B_R))$ (but not uniformly) and \eqref{timeBOUND}. Due to $f^\ep\in C([0,T];H^4(\R))$, in \eqref{regularizedM} it is easy to bound the linear terms. The nonlinear term can be written as
\begin{align*}
N(f)=&-C_{\ep}\Lambda^{1-\ep}f^\ep\\
&+\frac{\rho}{2\pi} PV\int_{\R^2} dy \frac{(\grad f(x)-\grad f(x-y))\cdot y}{|y|^{2-\ep}} ~ \frac{1-[1+(\Delta_yf(x))^2]^{3/2}}{[1+(\Delta_yf(x))^2]^{3/2}},
\end{align*}
and therefore
$$
|N(x,t)|\leq C(\ep)\|f^\ep\|_{H^4}(t),
$$
by Sobolev embedding.

The norm of $f_t^\ep\in W_*^{-2,\infty}(B_R)$ is given by
$$
\|f^\ep_t\|_{W_*^{-2,\infty}(B_R)}(t)
 = \sup_{\phi \in W_0^{2,1}(B_R): \| \phi \|_{W^{2,1}} \le 1} \left| \int_{\R}  dx ~ f_t^\ep(x,t) \phi(x) \right|,
$$
since $\phi$ vanishes on the boundary of $B_R$. Then we have
$$
I=\int_{B_R} \Lambda^{1-\ep}f(x) \phi(x)dx=\int_{\R^2} \Lambda^{1-\ep}f(x) \phi(x)dx=\int_{\R^2} f(x) \Lambda^{1-\ep}\phi(x)dx,
$$
and therefore
$$
|I|\leq \|f\|_{L^\infty}(t)\|\Lambda^{1-\ep}\phi\|_{L^1}.
$$
We split
$$
\Lambda^{1-\ep}\phi(x)=c\int_{\R^2}\frac{\phi(x)\!-\!\phi(x\!-\!y)}{|y|^{3-\ep}}dy=\int_{|y|>1}\!\!dy+\int_{|y|<1}\!\!dy=J_1(x)+J_2(x),
$$
so that
$$
\int_{\R^2}|J_1(x)|dx\leq \int_{|y|>1}\frac{dy}{|y|^{3-\ep}}\int_{\R^2} dx(|\phi(x)|+|\phi(x-y)|)\leq C\|\phi\|_{L^1(B_R)}.
$$
We rewrite $J_2$ as follows
$$
J_2(x)=c\int_{|y|<1}\frac{\phi(x)-\phi(x-y)-\grad\phi(x)\cdot y}{|y|^{3-\ep}}dy.
$$
We also consider the following identities
\begin{gather*}
\phi(x)-\phi(x-y)=\int_0^1 \grad\phi(x+(s-1)y)\cdot yds,
\\
\phi(x)\!-\!\phi(x\!-\!y)\!-\!\grad\phi(x)\cdot y=\int_0^1(s\!-\!1)ds\int_0^1dr\, y\cdot(\grad^2\phi(x+r(s\!-\!1)\al)\cdot y),
\end{gather*}
The expression for $J_2$ and these identities together yield
\begin{align*}
\int_{\R^2}|J_2(x)|dx&\leq \int_{|y|<1}|y|^{\ep-1}\int_0^1ds\int_0^1dr\int_{\R^2}dx ~ |\grad^2\phi(x+r(s-1)y)|\\
&\leq C\|\grad^2\phi\|_{L^1(B_R)}.
\end{align*}
We obtain
$$
\left\| \Lambda^{1-\ep}f^\ep \right\|_{W_*^{-2,\infty}(B_L)}
+
\left\|\Delta f^\ep\right\|_{W_*^{-2,\infty}(B_L)}
\le
C \left\| f^\ep \right\|_{L^{\infty}(\R^2)}
\le
C \left\| f_0 \right\|_{L^{\infty}(\R^2)}.
$$
For the last term in \eqref{regularizedM} we integrate by parts
$$
\int_{\R}  dx ~  \grad \phi(x) ~\cdot  PV\frac{\rho}{2\pi} PV\int_{\R^2} dy\frac{y}{|y|^{2-\ep}} ~ \frac{\Delta_yf(x)}{[1+(\Delta_yf(x))^2]^{1/2}}
$$
to realize that the splitting from \eqref{split.integral}  with $L=\delta=1$ allows us to conclude that the integral above is bounded by $C\| \phi \|_{W^{1,1}} \left\| f_0 \right\|_{L^{\infty}(\R)}$.

\section{Global existence for initial data in critical spaces}

This section is devoted to show global existence results for strong solutions of the Muskat contour equation in critical spaces.

\begin{thm}
Suppose that $f_0\in L^2$ and $\|f_0\|_{1}< k_0$ ($\|f_0\|_{1}< c_0$ for the 2D case).
Then there is a unique solution $f$ of Muskat with initial data $f_0$ that satisfies
$$
\|f\|_{L^2}(t)\leq \|f_0\|_{L^2},\quad \|f\|_1(t)+\mu\int_0^t ds\|f\|_2(s)\leq \|f_0\|_1,
$$
for $\mu>0$, a.e. $t\in [0,T]$ and any $T>0$. The time derivative of $f$ satisfies
$$
\|f_t\|_0(t)\leq C,\quad \int_0^Tds\|f_t\|_1(s)\leq C,
$$
where $C=C(\|f_0\|_1)$.
\end{thm}
\begin{rem}
The scale invariance for Muskat solutions $f^\lambda(x,t)=\frac1\lambda f(\lambda x,\lambda t)$ makes the following norms critical:
$$
\esssup_{t\in[0,T]}\|f\|_1(t),\quad \int_0^T ds\|f\|_2(s).
$$
The control of these norms gives in particular solutions such that
$$
\esssup_{t\in[0,T]}\|\grad f\|_{C_0}(t)+\mu\int_0^T ds\|\grad^2f\|_{C_0}(s)\leq \|f_0\|_1,
$$
where $C_0$ is the space of continuous functions vanishing at infinity.
\end{rem}
Proof: For $f_0$ such that $\|f_0\|_{1}< k_0$ we proceed as before to obtain the following a priori bound
$$
\frac{d}{dt}\|f\|_{1}(t)\leq 0,\quad \|f\|_{1}(t)<k_0.
$$
Due to
$$
1>\pi\sum_{n\geq 1} (2n+1)a_n\|f_0\|_1^{2n}=1-\mu\geq \pi\sum_{n\geq 1} (2n+1)a_n\|f\|_1^{2n}(t)
$$
for $0<\mu<1$, we find
$$
\frac{d}{dt}\|f\|_{1}(t)\leq -\mu\|f\|_{2}(t),
$$
(see Section 3 for details) and time integration gives the desired a priori bound
$$
\|f\|_1(t)+\mu\int_0^t\|f\|_2(s)ds\leq \|f_0\|_1.
$$
We also find
\begin{align}
\begin{split}\label{ftapriori}
\int d\xi |\hat{f}_t(\xi)|&\leq\! \int d\xi|\xi||\widehat{f}(\xi)|\!+\!\int d\xi |\mathcal{F}(N(f))(\xi)|\leq \|f\|_1(1\!+\!\pi\D\sum_{n\geq 1} a_n\|f\|_1^{2n})\leq C(\|f_0\|_1),
\end{split}
\end{align}
and similarly
\begin{align}
\begin{split}\label{ftapriori2}
\int_0^T\!\!\!\int dt d\xi |\xi||\hat{f}_t(\xi)|&\leq \int_0^T\!\!\!\int dtd\xi \Big(|\xi|^2|\widehat{f}(\xi)|+|\xi||\mathcal{F}(N(f))(\xi)|\Big)\\
&\leq \int_0^Tdt\|f\|_2(t)(1\!+\!\pi\D\sum_{n\geq 1}(2n\!+\!1) a_n\|f\|_1^{2n})\leq  C(\|f_0\|_1).
\end{split}
\end{align}

Next we would like to find a bona fide solution of Muskat satisfying those bounds.  We consider the following regularized model
$$
f_t^\ep=\zeta_\ep*(T(f^\ep)),\quad f^\ep(x,0)=(\zeta_\ep*f_0)(x),
$$
where
$$
T(f^\ep)(x)=\frac{1}{2\pi}\grad_x\cdot PV\int_{\R^2} dy\frac{y}{|y|^{2}} ~ \frac{\Delta_y(\zeta_\ep*f^\ep)(x)}{[1+(\Delta_y(\zeta_\ep*f^\ep)(x))^2]^{1/2}},
$$
(we take $\rho=1$ for the sake of simplicity and $\zeta_\ep$ given by \eqref{zeta}). Local existence can be shown as in \cite{DP} for regular initial data, since $f_0\in L^2$ it is easy to find $\zeta_\ep*f_0\in H^k$ for any $k\geq 0$. Then, as in Section 2, it is possible to obtain an $L^2$ maximum principle:
$$
\frac{1}{2}\frac{d}{dt}\|f^\ep\|^2_{L^2}=-\int_{\R^2}\int_{\R^2}
\frac{1}{|y|} \Big(1-\frac{1}{[1+(\Delta_{y}(\zeta_\ep*f^\ep(x)))^2]^{1/2}}\Big)dx dy.
$$
Due to $\|\zeta_\ep*f^\ep\|_{C^{2,\delta}}\leq \|\zeta_\ep*f^\ep\|_{H^4}\leq C(\ep)\|f^\ep\|_{L^2}\leq C(\ep)\|f_0\|_{L^2}$ it is possible to get global in time bounds and therefore global existence for $f^\ep\in C([0,T];H^k)$ for any $k\geq 3$ and any $T>0$ (see Section 2). Proceeding as before we find
$$
\|f^\ep\|_1(t)+\mu\int_0^t ds\|\zeta_\ep*f^\ep\|_2(s)\leq \|f_0\|_1.
$$

Next, we will take the limit as $\ep\to 0$. We will find strong and weak limits so most of the time the argument will be up to various subsequences. All of them will be denoted by $f^{\ep_n}$ by abuse of notation.

In particular $f^\ep$ is uniformly bounded in $L^\infty([0,T];L^2)$ so that there exists a subsequence $\{f^{\ep_n}\}$ which converges in the weak* topology of $L^\infty([0,T];L^2)$ to $f$. The subsequence $\{\widehat{f}^{\ep_n}\}$ is also uniformly bounded in $L^2([0,T]\times \R^2)$ so there exists a subsequence $\{\widehat{f}^{\ep_n}\}$ that converges weakly to $\widehat{f}\in L^2([0,T]\times \R^2)$. Then it is easy to check that $(\zeta_{\ep_n} * f^{\ep_n})\widehat{\,\,}(\xi,t)=\widehat{\zeta}(\ep_n\xi)\widehat{f}^{\ep_n}(\xi,t)$ converges weakly to $\widehat{f}\in L^2([0,T]\times \R^2)$.

We use Mazur's lemma to conclude that a convex combination
$$
G_n(\xi,t)=(G_n^1(\xi,t),G_n^2(\xi,t))=\sum_{k=n}^{N(n)}\lambda_k(\widehat{f}^{\ep_k}(\xi,t),\widehat{\zeta}(\ep_k\xi)\widehat{f}^{\ep_k}(\xi,t)),
$$
of $(\widehat{f}^{\ep_n}(\xi,t),\widehat{\zeta}(\ep_n\xi)\widehat{f}^{\ep_n}(\xi,t))$ with $(\cdot,\cdot)$ denoting a vector and
$$
\lambda_k\geq 0, \quad \sum_{k=n}^{N(n)}\lambda_k=1,
$$
converges strongly to $(\widehat{f},\widehat{f})$ in $(L^2([0,T]\times \R^2))^2$. We extract a subsequence (denoted by $G_n$) to get that $G_n(\xi,t)$ converges to $(\widehat{f}(\xi,t),\widehat{f}(\xi,t))$ pointwise for almost every $(\xi,t)\in \R^2\times [0,T]$. Therefore for $t\in [0,T]\smallsetminus \Omega$ with $|\Omega|=0$ we find that $G_n^1(\xi,t)$ converges to $\widehat{f}(\xi,t)$ pointwise for almost every $\xi\in\R^2$. We use Fatou's lemma to conclude that for $t\in [0,T]\smallsetminus \Omega$ and
$$
M(t)=\|f\|_1(t)+\mu\int_0^t ds\|f\|_2(s),
$$
the following holds
\begin{align*}
M(t)&\leq \liminf_{n\to\infty}\Big(\int d\xi|\xi||G^1_n(\xi,t)|+\mu \int_0^t ds \int d\xi|\xi|^2|G^2_n(\xi,s)|\Big)\\
    &\leq \liminf_{n\to\infty}\sum_{k=n}^{N(n)}\lambda_k\Big(\int \!d\xi|\xi| |\widehat{f}^{\ep_k}(\xi,t)|+\mu \int_0^t\!ds\!\int\! d\xi|\xi|^2 |\widehat{\zeta}(\ep_k\xi)\widehat{f}^{\ep_k}(\xi,s)|\Big)\\
 &\leq \liminf_{n\to\infty}\sum_{k=n}^{N(n)}\lambda_k \Big(\|f^{\ep_k}\|_1(t)+\mu\int_0^t ds\|\zeta_{\ep_k}*f^{\ep_k}\|_2(s)\Big)\leq \|f_0\|_1.
\end{align*}
Therefore
$$
\esssup_{t\in[0,T]} \|f\|_{1}(t)+\mu\int_0^T ds\|f\|_2(s)\leq \|f_0\|_1.
$$

In order to find that the limit function $f$ satisfies Muskat equation we claim that $f$ is a weak solution. Then the regularity of $f$ allows to conclude that it is in fact a strong solution. We will follow the arguments in Section 5 and Lemma 5.3 to get strong convergence in $L^\infty$. We just need to bound $f_t^{\ep_n}$ uniformly in $L^\infty([0,T];W_*^{-2,\infty}(B_R))$. But
$$
\|f_t^{\ep_n}\|_{W_*^{-2,\infty}(B_R)}(t)\leq \|f_t^{\ep_n}\|_{L^\infty(B_R)}(t)\leq \|f_t^{\ep_n}\|_{0}(t)\leq C(\|f_0\|_1),
$$
since the last inequality can be obtained as we did in the a priori bound \eqref{ftapriori}.  Since $\{f^{\ep_n}\}$ satisfies
\begin{multline}\label{weaksolep}
\int_0^T\!\!\int_{\R^2}\eta_t(x,t)f^{\ep_n}(x,t)dxdt+\int_{\R^2}\eta(x,0)(\zeta_{\ep_n}*f_0)(x)dx
\\
=\int_0^T\!\!\int_{\R^2}\!\!\!\!\grad_x(\zeta_{\ep_n}\!*\!\eta)(x,t)\cdot\frac{1}{2\pi}PV\!\!\int_{\R^2}\!\frac{y}{|y|^2}
\frac{\Delta_y (\zeta_{\ep_n}\!*\!f^{\ep_n})(x,t)dy dxdt}{[1+(\Delta_y (\zeta_{\ep_n}\!*\!f^{\ep_n})(x,t))^2]^{1/2}},
\end{multline}
we can pass to the limit as $\ep_n\to 0$ and the strong convergence gives $f$ as a weak Muskat solution.

Now we have $f$ a strong Muskat solution due to its regularity and we can find bounds \eqref{ftapriori} and \eqref{ftapriori2}. In order to end the result we just need to get uniqueness.

We consider two Muskat solutions $f_1$ and $f_2$ with the above properties and $f_1(x,0)=f_2(x,0)=f_0(x)$. Then for the difference $f=f_1-f_2$ we find
$$
\frac12\frac{d}{dt}\|f\|_{L^2}^2(t)=I+II+III,
$$
where
$$
I=\frac{1}{2\pi}\int f(x)\grad f(x)\cdot PV\int\frac{y}{[|y|^2+(f_1(x)-f_1(x-y))^2]^{3/2}}dydx,
$$
$$
II=\frac{-1}{2\pi}\int f(x)PV\int\frac{\grad
f(x-y)\cdot y}{[|y|^2+(f_1(x)-f_1(x-y))^2]^{3/2}}dydx,
$$
and
$$
III=\frac{1}{2\pi}\int f(x)\int \grad \Delta_y f_2(x)\cdot \frac{y}{|y|^2}
\Big(\frac{1}{[1+(\Delta_y f_1(x))^2]^{3/2}}-\frac{1}{[1+(\Delta_y f_2(x))^2]^{3/2}}\Big)dydx.
$$

We integrate by parts in $I$ to get
$$
I=3\int |f(x)|^2 A(x)dx,\quad \mbox{for}\quad A(x)=\frac{1}{4\pi}PV\int \frac{y}{|y|^3}\cdot\frac{\grad_x\Delta_y f_1(x)\Delta_yf_1(x)}{[1+(\Delta_yf_1(x))^2]^{5/2}}dy.
$$
Next, we bound as follows
$$
|A(x)|\leq \int d\xi|\mathcal{F}(A)(\xi)|.
$$
In order to deal with $\mathcal{F}(A)(\xi)$ we proceed as for $N(f)$ in \eqref{muskatEQ2d}. Since $z(1+z^2)^{-5/2}=\sum_{n\geq0}b_nz^{2n+1}$ for $|z|<1$ we can obtain
$$
\int d\xi|\mathcal{F}(A)(\xi)|\leq \frac{\pi}{4}\|f_1\|_2(t)\sum_{n\geq0}|b_n|(\|f\|_1(t))^{2n+1}\leq C(\|f_0\|_1)\|f_1\|_2(t).
$$
This yields
$$
I\leq  C(\|f_0\|_1)\|f_1\|_2(t)\|f\|_{L^2}^2(t).
$$

In the term $II$ we write $\grad f(x-y)=\grad_y(f(x)-f(x-y))$ and integrate by parts in $y$ to find $II=II_1+II_2$ where
$$
II_1=\frac{-1}{2\pi}\int f(x)PV\int\frac{(f(x)-f(x-y))}{[|y|^2+(f_1(x)-f_1(x-y))^2]^{3/2}}dydx,
$$
and
\begin{multline*}
II_2=\frac{3}{2\pi}\int f(x)PV\int (f(x)-f(x-y))\\
\times\frac{(f_1(x)-f_1(x-y))(f_1(x)-f_1(x-y)-\grad f_1(x-y)\cdot y)}{[|y|^2+(f_1(x)-f_1(x-y))^2]^{5/2}}dydx.
\end{multline*}
One could symmetrize $II_1$ to get
$$
II_1=\frac{-1}{4\pi}\int \int\frac{(f(x)-f(y))^2}{[|x-y|^2+(f_1(x)-f_1(y))^2]^{3/2}}dydx\leq 0.
$$
For $II_2$ we split further $II_2=II_2^1+II_2^2$ where
$$
II^1_2=3\int |f(x)|^2 B(x)dx,\quad \mbox{for}\quad
B(x)=\frac{1}{2\pi}PV\int B(x,y)dy,
$$
$$
B(x,y)=\frac{1}{|y|^3}\frac{\Delta_yf_1(x)(\Delta_yf_1(x)-\grad f_1(x-y)\cdot \frac{y}{|y|})}{[1+(\Delta_yf_1(x))^2]^{5/2}},
$$
and
$$
II^2_2=\frac{-3}{2\pi}\int f(x) PV\int f(x-y)B(x,y)dydx.
$$
Since
$$
f_1(x)-f_1(x-y)=\int_0^1\grad f_1(x+(s-1)y) ds \cdot y,
$$
we denote
$$
\grad_x\Delta_y^{s}  f_1(x)=\frac{\int_0^1\grad f_1(x+(s-1)y) ds-\grad f_1(x-y)}{|y|}
$$
to rewrite $B$ as follows
$$
B(x)=\frac{1}{2\pi}PV\int \frac{y}{|y|^3}\cdot\frac{\grad_x\Delta_y^s f_1(x)\Delta_yf_1(x)}{[1+(\Delta_yf_1(x))^2]^{5/2}}dy.
$$
At this point it is easy to find that $B$ and $A$ are similar in such a way that an analogous analysis allows us to obtain
$$
|B(x)|\leq\frac\pi2\|f_1\|_2(t)\sum_{n\geq0}|b_n|(\|f_1\|_1(t))^{2n+1}\leq C(\|f_0\|_1)\|f_1\|_2(t).
$$

It is possible to symmetrize $II_2^2$ as follows
$$
II^2_2=\frac{-3}{4\pi}\int f(x) PV\int [f(x-y)B(x,y)+f(x+y)B(x,-y)]dydx,
$$
and to use Parseval's identity in order to obtain
\begin{multline*}
II^2_2
=
\int d\xi \overline{\widehat{f}(\xi)}\D\sum_{n\geq 0}b_n\int d\xi_1\widehat{f}(\xi-\xi_1)\int d\xi_2\cdots\int d\xi_{2n+2}
\\
\times(\xi_1\!-\!\xi_2)\Big(\prod_{j=1}^{2n+1}\hat{f_1}(\xi_j\!-\!\xi_{j+1})\Big)
\hat{f_1}(\xi_{2n+2})\cdot J_n.
\end{multline*}
The integral $J_n=J_n(\xi,\xi_1,\ldots,\xi_{2n+2})$ reads
$$
J_n=\frac{-3i}{4\pi}PV\int \frac{y}{|y|^3}(M_n(y)-M_n(-y))dy,
$$
where $M_n(y)=M_n(\xi,\xi_1,\ldots,\xi_{2n+2},y)$ is given by
\begin{align*}
M_n(y)=e^{-i(\xi-\xi_1)\cdot y}&
e^{-i(\xi_1-\xi_2)\cdot y}\frac{e^{is(\xi_1-\xi_2)\cdot y}-1}{|y|} m(\xi_2\!-\!\xi_{3},y)\\
&\times m(\xi_3\!-\!\xi_{4},y)\ldots m(\xi_{2n+1}\!-\!\xi_{2n+2},y) m(\xi_{2n+2},y),
\end{align*}
using the operator $\grad_x\Delta_y^s$ in the $B(x,y)$ formula. With the PV cancelation we get
$$
|J_n|\leq \frac{3\pi}{2}\prod_{j=1}^{2n+1}|\xi_j-\xi_{j+1}||\xi_{2n+2}|,
$$
and therefore
\begin{align*}
II_2^2&\leq \frac{3\pi}{2}\int d\xi |\widehat{f}(\xi)| \D\sum_{n\geq 0}|b_n| \big[|\hat{f}|\underbrace{*(|\cdot|^2|\hat{f_1}|)*(|\cdot||\hat{f_1}|)\ldots*}_{2n+2 \mbox{ \footnotesize
 convolutions}}(|\cdot||\hat{f_1}|)\big](\xi)\\
 &\leq \frac{3\pi}{2}\|f\|_{L^2}^2\|f_1\|_{2}\D\sum_{n\geq 0}|b_n| \|f_1\|_1^{2n+1}\leq C(\|f_0\|_1)\|f_1\|_2(t)\|f\|_{L^2}^2(t).
\end{align*}
Above we use Schwarz's and Young's inequalities. It yields the desired estimate for $II$:
$$
II\leq  C(\|f_0\|_1)\|f_1\|_2(t)\|f\|_{L^2}^2(t).
$$

 We expand $III$ to obtain
$$
III=\frac{1}{2\pi}\int f(x)\sum_{n\geq 1}(-1)^na_n\int \frac y{|y|^2}\cdot\grad_x\Delta_y f_2(x)
[(\Delta_y f_1(x))^{2n}-(\Delta_y f_2(x))^{2n}]dydx.
$$
Since
\begin{multline*}
III=\frac{1}{2\pi}\int f(x)\sum_{n\geq 1}(-1)^na_n\int \frac y{|y|^2}\cdot\grad_x\Delta_y f_2(x)\Delta_y f(x)\\
\times\sum_{j=1}^{2n}(\Delta_y f_1(x))^{2n-j}(\Delta_y f_2(x))^{j-1}dydx,
\end{multline*}
we split further $III=III_1+III_2$ to find
\begin{multline*}
III_1=\frac{1}{2\pi}\int |f(x)|^2\sum_{n\geq 1}(-1)^na_nPV\int \frac y{|y|^3}\cdot\grad_x\Delta_y f_2(x)\\
\times\sum_{j=1}^{2n}(\Delta_y f_1(x))^{2n-j}(\Delta_y f_2(x))^{j-1}dydx,
\end{multline*}
and
\begin{multline*}
III_2=\frac{1}{2\pi}\int f(x)\sum_{n\geq 1}(-1)^na_n PV\int f(x-y)\frac y{|y|^3}\cdot\grad_x\Delta_y f_2(x)\\
\times\sum_{j=1}^{2n}(\Delta_y f_1(x))^{2n-j}(\Delta_y f_2(x))^{j-1}dydx.
\end{multline*}
We can proceed as before to get
$$
III_1\leq \|f\|_{L^2}^2(t)\|f_2\|_2(t)\frac\pi2\sum_{n\geq 1}2n a_n\|f_0\|_1^{2n-1}\leq C(\|f_0\|_1)\|f_2\|_2(t)\|f\|_{L^2}^2(t).
$$
We deal with $III_2$ as with $II_2^2$:
$$
III_2\leq \|f\|_{L^2}^2(t)\|f_2\|_2(t)\frac\pi2\sum_{n\geq 1}2n a_n\|f_0\|_1^{2n-1}\leq C(\|f_0\|_1)\|f_2\|_2(t)\|f\|_{L^2}^2(t).
$$
We obtain finally
$$
\frac12\frac{d}{dt}\|f\|_{L^2}^2(t)\leq C(\|f_0\|_1)(\|f_1\|_2(t)+\|f_2\|_2(t))\|f\|_{L^2}^2(t),
$$
and time integration provides
$$
\|f\|_{L^2}^2(t)\leq \|f_0\|_{L^2}^2\exp\Big(\D C(\|f_0\|_1)\int_0^t(\|f_1\|_2(s)+\|f_2\|_2(s))ds\Big),
$$
to find $f=0$.  This completes the proof of uniqueness.
\hfill$\Box$

Next we consider the following norms for $s,\,p\geq 1$ given by
$$
\|f\|^p_{s,p}=\int_{\R}|\xi|^{sp}|\widehat{f}(\xi)|^p d\xi,
$$
with the homogeneous space
$$
\mathcal{F}^{s,p}=\{T\in\mathcal{S}'(\R): \mbox{$\widehat{T}$ is a function and }\|T\|_{s,p}<\infty \}.
$$
We provide the following result:

\begin{thm}
Suppose that $f_0\in L^2\cap\mathcal{F}^{1,1}\cap\mathcal{F}^{2-\frac1p,p}$ with $p>1$ and
$$
\|f_0\|_{1}=\|f_0\|_{1,1}< k_0,\quad (\|f_0\|_1<c_0\mbox{ in 2D}).
$$
Then there is a unique solution $f$ of Muskat with initial data $f_0$ that satisfies
$$
f\in L^{\infty}([0,T];L^2\cap \mathcal{F}^{1,1}\cap \mathcal{F}^{2-\frac1p,p})\cap L^{1}([0,T];\mathcal{F}^{2,1})\cap L^{p}([0,T];\mathcal{F}^{2,p})
$$
for any $T>0$. The time derivative of $f$ satisfies
$$
f_t\in L^{\infty}([0,T];\mathcal{F}^{0,1}\cap\mathcal{F}^{1-\frac1p,p})\cap L^{1}([0,T];\mathcal{F}^{1,1})\cap L^{p}([0,T];\mathcal{F}^{1,p}).
$$
\end{thm}

\begin{rem}
We would like to point out that all the homogeneous norms $$\sup_{[0,T]}\|f\|_{2-\frac1p,p},\quad\mbox{and}\quad \Big(\int_0^T\|f\|_{2,p}^p(t)dt\Big)^{1/p}$$ are critical in 2D under the scale invariant for Muskat contour equation. In 3D these norms are supercritical due to the fact that, for example,
$$
\sup_{[0,T]}\|f\|_{2,2}(t)
$$
is critical.

In particular, we give a new result in Sobolev spaces taking $p=2$. In the case $1<p<2$ it is possible to use Hausdorff-Young inequality to find
$f\in L^{\infty}([0,T],W^{2-\frac1p,\frac p{p-1}})\cap L^{p}([0,T];W^{2,\frac p{p-1}})$

 \end{rem}

Proof: For $f_0$ such that $\|f_0\|_{1}< k_0$ we proceed as before to obtain a priori estimates. Next we check the evolution of
\begin{align*}
\frac1p\frac{d}{dt}\|f\|_{2-\frac1p,p}^p&=\int_{\R} d\xi |\xi|^{2p-1}|\widehat{f}(\xi)|^{p-1}(\widehat{f}_t(\xi)\overline{\widehat{f}(\xi)}+\widehat{f}(\xi)\overline{\widehat{f}_t(\xi)})/(2|\widehat{f}(\xi)|)\\
&\leq -\int_{\R}d\xi|\xi|^{2p}|\widehat{f}(\xi)|^{p}+\int_{\R}d\xi|\xi|^{2p-1}|\widehat{f}(\xi)|^{p-1}|\mathcal{F}(N(f))(\xi)|\\
&=-\int_{\R}d\xi|\xi|^{2p}|\widehat{f}(\xi)|^{p}+I
\end{align*}
We bound as follows
\begin{align*}
I&\leq \int_{\R} d\xi |\xi|^{2(p-1)}|\widehat{f}(\xi)|^{p-1}|\xi|\pi\D\sum_{n\geq 1}a_n\int_{\R^2} d\xi_1\cdots\int_{\R^2} d\xi_{2n} ~
\\
&\qquad\quad\times  |\xi-\xi_1||\hat{f}(\xi-\xi_{1})| \prod_{j=1}^{2n-1}  |\xi_j -\!\xi_{j+1}||\hat{f}(\xi_j -\!\xi_{j+1})||\xi_{2n}||\hat{f}(\xi_{2n})|.
\end{align*}
The inequality $|\xi|\leq |\xi-\xi_1|+|\xi_1-\xi_2|+...+|\xi_{2n-1}-\xi_{2n}|+|\xi_{2n}|$ gives
\begin{align*}
I\leq \int_{\R} d\xi |\xi|^{2(p-1)}|\widehat{f}(\xi)|^{p-1}\pi\D\sum_{n\geq 1}(2n+1)a_n \big[(|\cdot|^2|\hat{f}|)\underbrace{*(|\cdot||\hat{f}|)*\ldots*}_{2n \mbox{ \footnotesize
 convolutions}}(|\cdot||\hat{f}|)\big](\xi).
\end{align*}
H\"older and Young's inequalities yield
\begin{align*}
I&\leq \|f\|_{2,p}^{p-1}\pi\D\sum_{n\geq 1}(2n+1)a_n \|(|\cdot|^2|\hat{f}|)\underbrace{*(|\cdot||\hat{f}|)*\ldots*}_{2n \mbox{ \footnotesize
 convolutions}}(|\cdot||\hat{f}|)\|_{L^p}\\
&\leq \|f\|_{2,p}^{p-1}\pi\D\sum_{n\geq 1}(2n+1)a_n \|f\|_{2,p}\|f\|_1^{2n}.
\end{align*}
Due to
$$
1>\pi\D\sum_{n\geq 1}(2n+1) a_n\|f_0\|_1^{2n}=1-\mu\geq \pi\D\sum_{n\geq 1}(2n+1)a_n \|f\|_1^{2n}(t)
$$
for $0<\mu<1$, we find
$$
\frac1p\frac{d}{dt}\|f\|_{2-\frac1p,p}^p(t)\leq -\mu\|f\|_{2,p}^p(t),
$$
and time integration gives the following a priori bound
$$
\|f\|_{2-\frac1p,p}^p(t)+p\mu\int_0^tds\|f\|_{2,p}^p(s)\leq \|f_0\|_{2-\frac1p,p}^p.
$$
We also find
\begin{align*}
\int d\xi |\hat{f}_t(\xi)|&\leq \int_{\R}d\xi|\xi||\widehat{f}(\xi)|+\int_{\R}d\xi |\mathcal{F}(N(f))(\xi)|\\
&\leq \|f\|_1(1+\pi\D\sum_{n\geq 1}a_n \|f\|_1^{2n})\leq 2.
\end{align*}
For $g\in L^{\frac p{p-1}}$ and $\|g\|_{L^{\frac p{p-1}}}\leq 1$ we find
\begin{align*}
\int_{\R}d\xi g(\xi) |\xi|^{1-\frac1p}\widehat{f}_t(\xi)&\leq \int_{\R}d\xi|g(\xi)|(|\xi|^{2-\frac1p}|\widehat{f}(\xi)|+|\xi|^{1-\frac1p}|\mathcal{F}(N(f))(\xi)|)\\
&\leq \|f\|_{2-\frac1p,p}+J.
\end{align*}
Using that $|\xi|^{1-\frac1p}\leq (|\xi-\xi_1|^{1-\frac1p}+|\xi_1-\xi_2|^{1-\frac1p}+...+|\xi_{2n-1}-\xi_{2n}|^{1-\frac1p}+|\xi_{2n}|^{1-\frac1p})$ we get
\begin{align*}
J&\leq \int_{\R} d\xi |g(\xi)| 2\D\sum_{n\geq 1}(2n+1) \big[(|\cdot|^{2-\frac1p}|\hat{f}|)\underbrace{*(|\cdot||\hat{f}|)*\ldots*}_{2n \mbox{ \footnotesize
 convolutions}}(|\cdot||\hat{f}|)\big](\xi)\\
 &\leq \pi\D\sum_{n\geq 1}(2n+1)^{1+\frac1p}a_n\|f\|_1^{2n}\|f\|_{2-\frac1p,p}\leq C(\|f_0\|_1)\|f\|_{2-\frac1p,p}.
\end{align*}
Therefore duality provides
$$
\|f_t\|_{1-\frac1p,p}\leq C(\|f_0\|_1)\|f_0\|_{2-\frac1p,p}.
$$
An analogous approach provides
$$
\int_0^Tdt\|f_t\|_{1,p}^p(t)\leq \frac{C(\|f_0\|_1)}{p\mu}\|f_0\|^p_{2-\frac1p,p}.
$$

In order to find bona fide solutions of the system we regularize the initial data as in the previous theorem. We then obtain the same  a priori bounds for the regularized solution as above. We pass to the limit to find a global-in-time solution. We find weak* and weak convergence of the regularized system to the Muskat solution in the weak* and weak topology of $L^{\infty}([0,T];\mathcal{F}^{2-\frac1p,p})$ and $L^{p}([0,T];\mathcal{F}^{2,p})$ respectively using that for $p>1$ the spaces $L^p$ are reflexive.

\subsection*{{\bf Acknowledgments}}

\smallskip

PC was partially supported by NSF
grants DMS-1209394 and DMS-1265132. DC and FG were partially supported by MCINN grant MTM2011-26696 (Spain).
FG acknowledges support from the
Ram\'on y Cajal program. LRP was partially supported by MCINN grant MTM2012-05622 (Spain). RMS was partially supported by NSF grants DMS-1200747, DMS-0901463, and an Alfred
P. Sloan Foundation Research Fellowship.

FG wish to thank Juan Casado-D\'iaz and David Ruiz for helpful discussions.

\newpage

\begin{quote}
\begin{tabular}{l}
\textbf{Peter Constantin }\\
{\small Department of Mathematics}\\
{\small Princeton University}\\
{\small 1102 Fine Hall, Washington Rd, Princeton, NJ 08544, USA}\\
{\small Email: const@math.princeton.edu}
\end{tabular}
\end{quote}

%\begin{quote}
%\begin{tabular}{ll}
%\textbf{Diego C\'ordoba} &  \textbf{Francisco Gancedo}\\
%{\small Instituto de Ciencias Matem\'aticas} & {\small Department of Mathematics}\\
%{\small Consejo Superior de Investigaciones Cient\'ificas} & {\small University of Chicago}\\
%{\small Serrano 123, 28006 Madrid, Spain} & {\small 5734 University Avenue, Chicago, IL 60637}\\
%{\small Email: dcg@icmat.es} & {\small Email: fgancedo@math.uchicago.edu}
%\end{tabular}
%\end{quote}

\begin{quote}
\begin{tabular}{ll}
\textbf{Diego C\'ordoba}\\
{\small Instituto de Ciencias Matem\'aticas}\\
{\small Consejo Superior de Investigaciones Cient\'ificas}\\
{\small C/ Nicolas Cabrera, 13-15, 28049 Madrid, Spain}\\
{\small Email: dcg@icmat.es}
\end{tabular}
\end{quote}

\begin{quote}
\begin{tabular}{ll}
\textbf{Francisco Gancedo}\\
{\small Departamento de An\'{a}lisis Matem\'{a}tico $\&$ IMUS}\\
{\small Universidad de Sevilla}\\
{\small C/ Tarfia s/n, Campus Reina Mercedes, 41012 Sevilla, Spain}\\
{\small Email: fgancedo@us.es}
\end{tabular}
\end{quote}

\begin{quote}
\begin{tabular}{ll}
\textbf{Luis Rodr\'iguez-Piazza}\\
{\small Departamento de An\'{a}lisis Matem\'{a}tico $\&$ IMUS}\\
{\small Universidad de Sevilla}\\
{\small C/ Tarfia s/n, Campus Reina Mercedes, 41012 Sevilla, Spain}\\
{\small Email: piazza@us.es}
\end{tabular}
\end{quote}

\begin{quote}
\begin{tabular}{ll}
\textbf{Robert M. Strain}\\
{\small Department of Mathematics}\\
{\small University of Pennsylvania}\\
{\small David Rittenhouse Lab}\\
{\small 209 South 33rd Street, Philadelphia, PA 19104, USA} \\ %-6395, USA}\\
{\small Email: strain@math.upenn.edu}
\end{tabular}
\end{quote}

\end{document}